\documentclass[aos,MSNbibl,seceqn,nameyear,dvips]{arximspdf}
\usepackage{dcolumn}
\usepackage{graphicx}

\doi{10.1214/11-AOS895}
\volume{39}
\issue{4}
\pubyear{2011}
\firstpage{2131}
\lastpage{2163}

\makeatletter

\newcolumntype{d}[1]{D{.}{.}{#1}}

\newtheorem{theorem}{Theorem}[section]
\newtheorem{lem}{Lemma}[section]
\newproclaim{asm}{Assumption}[section]
\newproclaim{rem}{Remark}[section]

\newcommand{\f}{\frac}
\newcommand{\p}{\partial}
\newcommand{\e}{\varepsilon}

\makeatother

\begin{document}
\begin{frontmatter}

\title{Global self-weighted and local quasi-maximum exponential likelihood estimators for ARMA--GARCH/IGARCH models}
\runtitle{QMELE for ARMA--GARCH/IGARCH models}

\begin{aug}
\author[A]{\fnms{Ke} \snm{Zhu}\ead[label=e1]{mazkxaa@gmail.com}}
and
\author[A]{\fnms{Shiqing} \snm{Ling}\corref{}\thanksref{t1}\ead[label=e2]{maling@ust.hk}}
\runauthor{K. Zhu and S. Ling}
\affiliation{Hong Kong University of Science and Technology}
\address[A]{Department of Mathematics\\
Hong Kong University\\
\quad of Science and Technology\\
Clear Water Bay, Kowloon\\
Hong Kong\\
\printead{e1}\\
\phantom{E-mail: }\printead*{e2}} %adresu isvedimo komanda gale!
\end{aug}

\thankstext{t1}{Supported in part by Hong Kong Research Grants Commission
Grants HKUST601607 and HKUST602609.}

% HISTORY:
\received{\smonth{1} \syear{2011}}

% ABSTRACT
%
\begin{abstract}
This paper investigates the asymptotic theory of the quasi-maxi\-mum
exponential likelihood estimators (QMELE) for ARMA--GARCH models. Under
only a fractional moment condition, the strong consistency and the
asymptotic normality of the global self-weighted QMELE are obtained.
Based on this self-weighted QMELE, the local QMELE is showed to be
asymptotically normal for the ARMA model with GARCH (finite variance)
and IGARCH errors. A formal comparison of two estimators is given for
some cases. A simulation study is carried out to assess the performance
of these estimators, and a real example on the world crude oil price is
given.
\end{abstract}

% KEYWORDS
%
\begin{keyword}[class=AMS]
\kwd{62F12}
\kwd{62M10}
\kwd{62P20}.
\end{keyword}
\begin{keyword}
\kwd{ARMA--GARCH/IGARCH model}
\kwd{asymptotic normality}
\kwd{global self-weighted/local
quasi-maximum exponential likelihood estimator}
\kwd{strong consistency}.
\end{keyword}

\end{frontmatter}

%s1 ###
%se1 #&#
\section{Introduction}

Assume that $\{y_t\dvtx t=0,\pm1,\pm2,\ldots\}$ is generated by the
ARMA--GARCH model
%
%e1.2 ###
%e1.1 ###
%e1.1 #&#
%e1.2 #&#
\begin{eqnarray}
\label{11}
y_t&=&\mu+\sum_{i=1}^p\phi_iy_{t-i}+\sum_{i=1}^q\psi_i\e_{t-i}+\e_t,
\\
\label{12}
\e_t&=&\eta_t\sqrt{h_t} \quad\mbox{and}\quad
h_t=\alpha_0+\sum_{i=1}^r\alpha_i\e_{t-i}^2+\sum_{i=1}^s\beta_ih_{t-i},
\end{eqnarray}
where $\alpha_0>0,\alpha_i\geq0$ $(i=1,\ldots,r),\beta_j\geq0$
$(j=1,\ldots,s)$, and $\eta_t$ is a sequence of i.i.d. random variables
with $E\eta_t=0$. As we all know, since \citet{r12} and \citet
{r5}, model (\ref{11})--(\ref{12}) has been widely used in economics and
finance; see \citet{r6}, \citet{r2},
\citet{r7} and \citet{r14}. The
asymptotic theory of the quasi-maximum likelihood estimator (QMLE)
was established by \citet{r26} and by \citet{r13} when $E\e
_t^4<\infty$. Under the strict stationarity
condition, the consistency and the asymptotic normality of the QMLE
were obtained by \citet{r21} and \citet{r29} for the
GARCH$(1,1)$ model, and by \citet{r4} and
\citet{r13} for the GARCH$(r, s)$ model. \citet{r15}
established the asymptotic theory of the QMLE for the GARCH
model when $E\e_t^{2}<\infty$, including both cases in which
$E\eta_t^{4}=\infty$ and $E\eta_t^{4}<\infty$. Under the geometric
ergodicity condition, \citet{r20} gave the asymptotic
properties of the modified QMLE for the first order AR--ARCH model.
Moreover, when $E|\e_t|^\iota<\infty$ for some $\iota>0$, the
asymptotic theory of the global self-weighted QMLE and the local
QMLE was established by \citet{r25} for model~\mbox{(\ref{11})--(\ref{12})}.

It is well known that the asymptotic normality of the QMLE requires
$E\eta_{t}^{4}<\infty$ and this property is lost when
$E\eta_{t}^{4}=\infty$; see \citet{r15}. Usually, the least
absolute deviation (LAD) approach can be used to reduce the moment
condition of $\eta_{t}$ and provide a robust estimator. The local
LAD estimator was studied by \citet{r31} and \citet{r22} for
the pure GARCH model, \citet{r9} for the double
AR(1) model, and \citet{r23} for the ARFIMA--GARCH model. The
global LAD estimator was studied by \citet{r16} for
the pure ARCH model and by \citet{r3} for the pure
GARCH model, and by \citet{r35} for the double AR($p$) model.
Except for the AR models studied by Davis, Knight and Liu (\citeyear{r11})
and \citet{r24} [see also
Knight (\citeyear{r18}, \citeyear{r19})], the nondifferentiable and
nonconvex objective function appears when one studies the LAD
estimator for the ARMA model with i.i.d. errors. By assuming the
existence of a $\sqrt{n}$-consistent estimator,
the asymptotic normality of the LAD estimator is established for the
ARMA model with i.i.d. errors
by \citet{r10} for the finite variance case and by \citet
{r30} for the infinite variance case;
see also Wu and Davis (\citeyear{r33}) for the noncausal or noninvertible
ARMA model. Recently, \citet{r36} proved the asymptotic
normality of the global LAD estimator for the finite/infinite variance
ARMA model with i.i.d. errors.

In this paper, we investigate the self-weighted quasi-maximum
exponential likelihood estimator (QMELE) for model (\ref{11})--(\ref{12}).
Under only a fractional moment condition of $\e_{t}$ with
$E\eta_{t}^{2}<\infty$, the strong consistency and the asymptotic
normality of the global self-weighted QMELE are obtained by using
the bracketing method in \citet{r32}. Based on this global
self-weighted QMELE, the local QMELE is showed to be asymptotically
normal for the ARMA--GARCH (finite variance) and --IGARCH models. A
formal comparison of two estimators is given for some cases.

%f1 ###
%fi1 #&#
\begin{figure}

\includegraphics{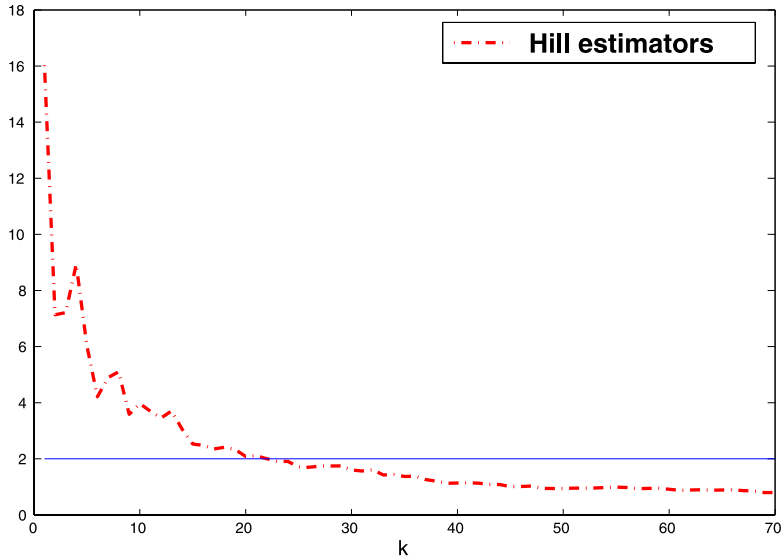}

\caption{The Hill estimators $\{\hat{\alpha}_{\eta}(k)\}$ for
$\{\hat{\eta_{t}}^{2}\}$.} \label{figure1}
\end{figure}

To motivate our estimation procedure, we revisit the GNP deflator
example of \citet{r5}, in which the GARCH model was proposed
for the first time. The model he specified is an AR(4)--GARCH$(1,1)$
model for the quarterly data from 1948.2 to 1983.4 with a total of
143 observations. We use this data set and his fitted model to
obtain the residuals $\{\hat{\eta}_{t}\}$. The tail index of
$\{\eta_{t}^{2}\}$ is estimated by Hill's estimator
$\hat{\alpha}_{\eta}(k)$ with the largest $k$ data of
$\{\hat{\eta}_{t}^{2}\}$, that is,
\[
\hat{\alpha}_{\eta}(k)=\frac{k}{\sum_{j=1}^{k} (\log\tilde{\eta
}_{143-j}-\log\tilde{\eta}_{143-k})},
\]
where $\tilde{\eta}_{j}$ is the $j$th order statistic of
$\{\hat{\eta}_{t}^{2}\}$. The plot of
$\{\hat{\alpha}_{\eta}(k)\}_{k=1}^{70}$ is given in Figure~\ref
{figure1}. From
this figure, we can see that $\hat{\alpha}_{\eta}(k)>2$ when $k\leq
20$, and $\hat{\alpha}_{\eta}(k)<2$ when $k>20$. Note that Hill's
estimator is not so reliable when~$k$ is too small. Thus, the tail
of $\{\eta_{t}^{2}\}$ is most likely less than 2, that is,
$E\eta_{t}^{4}=\infty$. Thus, the setup that $\eta_{t}$ has a finite
forth moment may not be suitable, and hence the standard QMLE
procedure may not be reliable in this case. The estimation procedure
in this paper only requires $E\eta_{t}^{2}<\infty$. It may provide a
more reliable alternative to practitioners. To further illustrate
this advantage, a simulation study is carried out to compare the
performance of our estimators and the self-weighted/local QMLE in
\citet{r25}, and a~new real example on the world crude oil price is
given in this paper.

This paper is organized as follows. Section \ref{sec2} gives our results on
the global self-weighted QMELE. Section \ref{sec3} proposes a local QMELE
estimator and gives its limiting distribution. The simulation
results are reported in Section~\ref{sec4}. A real example is given in
Section \ref{sec5}. The proofs of two technical lemmas are provided in
Section \ref{sec6}. Concluding remarks are offered in Section~\ref{sec7}. The
remaining proofs are given in the \hyperref[app]{Appendix}.

%s2 ###
%se2 #&#
\section{Global self-weighted QMELE}\label{sec2}
Let $\theta=(\gamma',\delta')'$ be the unknown parameter of model
(\ref{11})--(\ref{12}) and its true value be $\theta_{0}$, where
$\gamma=(\mu,\phi_1,\ldots,\phi_p,\allowbreak\psi_1,\ldots,\psi_q)'$ and
$\delta=(\alpha_0,\ldots,\alpha_r,\beta_1,\ldots,\beta_s)'$. Given the
observations $\{y_n,\ldots,\allowbreak y_1\}$ and the initial values
$Y_{0}\equiv\{y_0,y_{-1},\ldots\}$, we can rewrite the
parametric
model \mbox{(\ref{11})--(\ref{12})} as
%
%e2.2 ###
%e2.1 ###
%e2.1 #&#
%e2.2 #&#
\begin{eqnarray}\quad
\label{21}
\e_t(\gamma)&=&y_t-\mu-\sum_{i=1}^p\phi_iy_{t-i}-\sum_{i=1}^q\psi_i\e
_{t-i}(\gamma),\\
\label{22}
\eta_t(\theta)&=&\e_t(\gamma)/\sqrt{h_t(\theta)}
\quad\mbox{and}\nonumber\\[-8pt]\\[-8pt]
h_t(\theta)&=&\alpha_0+\sum_{i=1}^r\alpha_i\e_{t-i}^2(\gamma)+\sum
_{i=1}^s\beta_ih_{t-i}(\theta).\nonumber
\end{eqnarray}
Here, $\eta_t(\theta_0)=\eta_t$, $\e_t(\gamma_0)=\e_t$ and
$h_t(\theta_0)=h_t$.
The parameter space is $\Theta=\Theta_{\gamma} \times\Theta_{\delta
}$, where
$\Theta_{\gamma}\subset R^{p+q+1}$, $\Theta_{\delta}\subset
R^{r+s+1}_{0}$, $R=(-\infty, \infty)$ and $R_{0}=[0, \infty)$.
Assume that $\Theta_{\gamma}$ and $\Theta_{\delta}$ are compact and
$\theta_{0}$ is an interior point in $ \Theta$. Denote $\alpha(z)=
\sum_{i=1}^{r}\alpha_{i}z^{i}$, $\beta(z)= 1- \sum_{i=1}^{s}
\beta_{i}z^{i}$, $\phi(z)= 1-\sum_{i=1}^{p}\phi_{i}z^{i}$ and
$\psi(z)=1+ \sum_{i=1}^{q}\psi_{i}z^{i}$. We introduce the
following assumptions:
%
%as2.1 #&#
\begin{asm}\label{asm21}
For each
$\theta\in\Theta$, $\phi(z)\neq0$ and $\psi(z)\neq0$ when
$|z|\leq1$, and $\phi(z)$ and $\psi(z)$ have no common root with
$\phi_p\neq0$ or $\psi_q\neq0$.
\end{asm}
%
%as2.2 #&#
\begin{asm}\label{asm22}
For each $\theta\in\Theta$,
$\alpha(z)$ and $\beta(z)$ have no common root,
$\alpha(1)\neq1,\alpha_r+\beta_s\neq0$ and $\sum_{i=1}^s\beta_i<1$.
\end{asm}
%
%as2.3 #&#
\begin{asm}\label{asm23}
$\eta_t^{2}$ has a nondegenerate distribution with
$E\eta_t^{2}<\infty$.
\end{asm}

Assumption \ref{asm21} implies the stationarity, invertibility and
identifiability of mod\-el~(\ref{11}), and Assumption \ref{asm22} is the
identifiability condition for mo\-del~(\ref{12}).~Assumption \ref{asm23} is
necessary to ensure that $\eta_t^{2}$ is not almost surely (a.s.)
a~constant. When $\eta_{t}$ follows the standard double exponential
distribution, the weighted log-likelihood function (ignoring a
constant) can be written as follows:
%
%e2.3 ###
%e2.3 #&#
\begin{equation}\label{23}
L_{sn}(\theta)=\f{1}{n}\sum_{t=1}^n w_tl_t(\theta) \quad\mbox{and}\quad
l_t(\theta)=\log\sqrt{h_t(\theta)}+\f{|\e_t(\gamma)|}{\sqrt{h_t(\theta)}},
\end{equation}
where $w_t=w(y_{t-1},y_{t-2},\ldots)$ and $w$ is a measurable, positive
and bounded function on $R^{Z_0}$ with $Z_0=\{0,1,2,\ldots\}$. We look\vadjust{\goodbreak}
for the minimizer,
$\hat{\theta}_{sn}=(\hat{\gamma}_{sn}',
\hat{\delta}_{sn}')'$, of $L_{sn}(\theta)$ on $\Theta$, that is,
\[
\hat{\theta}_{sn}=\mathop{\arg\min}_{\Theta}L_{sn}(\theta).
\]
Since the weight $w_{t}$ only depends on $\{y_{t}\}$ itself and we
do not assume that~$\eta_{t}$ follows the standard double exponential distribution,
$\hat{\theta}_{sn}$ is called the self-weighted quasi-maximum
exponential likelihood estimator (QMELE) of $\theta_{0}$. When
$h_{t}$ is a constant, the self-weighted QMELE reduces to the
weighted LAD estimator of the ARMA model in \citet{r30} and
\citet{r36}.

The weight $w_{t}$ is to reduce the moment condition of $\e_t$ [see
more discussions in \citet{r25}], and it satisfies the following
assumption:
%
%as2.4 #&#
\begin{asm} \label{asm24}
$E[(w_t+w_t^2)\xi_{\rho t-1}^3]<\infty$ for any
$\rho\in(0,1)$, where $\xi_{\rho
t}=1+\sum_{i=0}^\infty\rho^i|y_{t-i}|$.
\end{asm}

\begin{figure}

\includegraphics{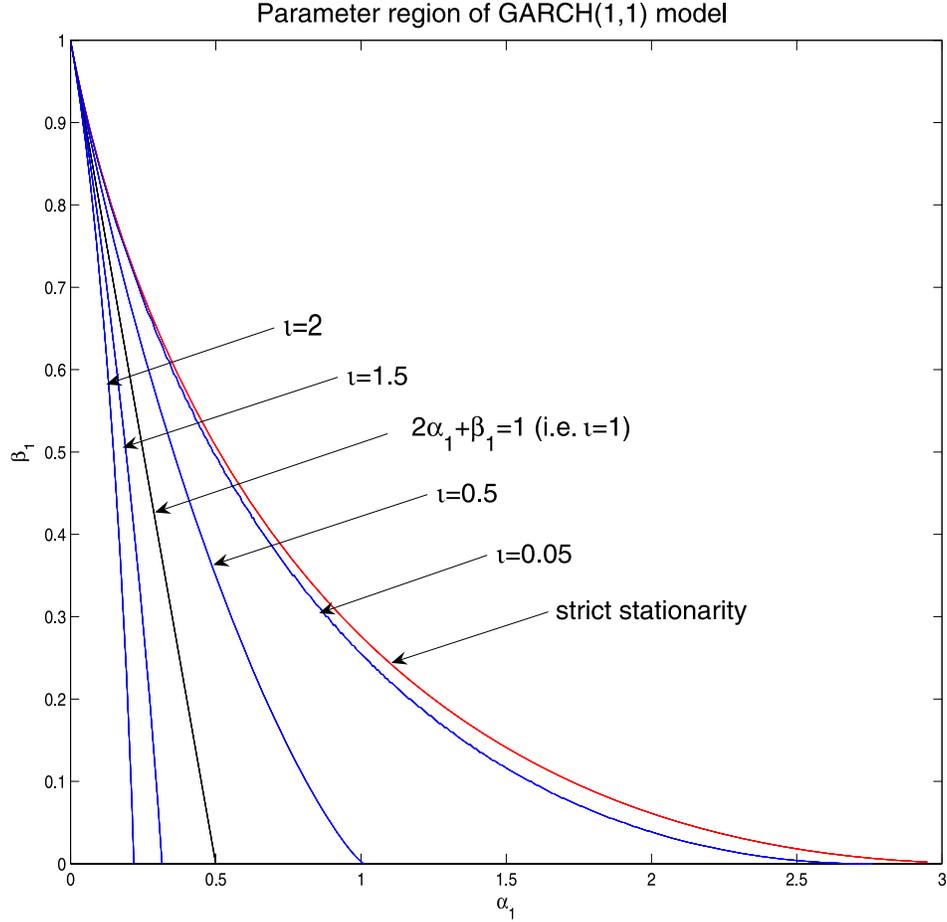}

\caption{The regions bounded by the indicated curves are for the
strict stationarity and for $E|\e_t|^{2\iota}<\infty$ with
$\iota=0.05, 0.5, 1, 1.5$ and $2$, respectively.} \label{figure2}
\end{figure}

When $w_{t}\equiv1$, the $\hat{\theta}_{sn}$ is the
global QMELE and it needs the moment condition
$E|\varepsilon_{t}|^{3}<\infty$ for its asymptotic normality, which
is weaker than the moment condition $E\varepsilon_{t}^4<\infty$ as
for the QMLE of $\theta_0$ in \citet{r13}. It is
well known that the higher is the moment condition of $\e_t$, the
smaller is the parameter space. Figure \ref{figure2} gives the strict
%
%f2 ###
%fi2 #&#
stationarity region and regions for $E|\e_t|^{2\iota}<\infty$ of the
GARCH$(1,1)$ model: $\e_t=\eta_t\sqrt{h_t}$ and
$h_t=\alpha_0+\alpha_1\e_{t-1}^{2}+\beta_1 h_{t-1}$, where
$\eta_t\sim \operatorname{Laplace}(0,1)$. From Figure \ref{figure2}, we can see that
the region
for $E|\e_t|^{0.1}<\infty$ is very close to the region for strict
stationarity of $\e_t$, and is much bigger than the region for
$E\e_t^{4}<\infty$.

Under Assumption \ref{asm24}, we only need a
fractional moment condition for the asymptotic property of
$\hat{\theta}_{sn}$ as follows:
%
%as2.5 #&#
\begin{asm} \label{asm25}
$E|\e_t|^{2\iota}<\infty$ for some $\iota>0$.
\end{asm}

The sufficient and necessary condition of Assumption \ref{asm25} is
given in
Theorem~2.1 of \citet{r25}. In practice, we can use Hill's estimator
to estimate the tail index of $\{y_{t}\}$ and its estimator may
provide some useful guidelines for the choice of $\iota$. For
instance, the quantity $2\iota$ can be any value less than the tail
index $\{y_{t}\}$. However, so far we do not know how to choose the
optimal $\iota$. As in \citet{r25} and \citet{r30}, we choose
the weight function $w_{t}$ according to $\iota$. When $\iota=1/2$
(i.e., $E|\e_{t}|<\infty$), we can choose the weight function as
%
%e2.4 ###
%e2.4 #&#
\begin{equation}\label{24}
w_{t}=\Biggl(\max\Biggl\{1,C^{-1}\sum_{k=1}^{\infty}
\frac{1}{k^{9}}|y_{t-k}|I\{|y_{t-k}|>C\}\Biggr\}\Biggr)^{-4},
\end{equation}
where $C>0$ is a constant. In practice, it works well when we select
$C$ as the 90\% quantile of data $\{y_{1},\ldots,y_{n}\}$. When\vadjust{\goodbreak}
$q=s=0$ (AR--ARCH model), for any $\iota>0$, the weight can be
selected as
\[
w_{t}=\Biggl(\max\Biggl\{1,C^{-1}\sum_{k=1}^{p+r}
\frac{1}{k^{9}}|y_{t-k}|I\{|y_{t-k}|>C\}\Biggr\}\Biggr)^{-4}.
\]
When $\iota\in(0,1/2)$ and $q>0$ or $s>0$, the weight function need
to be modified as follows:
\[
w_{t}=\Biggl(\max\Biggl\{1,C^{-1}\sum_{k=1}^{\infty}
\frac{1}{k^{1+8/\iota}}|y_{t-k}|I\{|y_{t-k}|>C\}\Biggr\}\Biggr)^{-4}.
\]
Obviously, these weight functions satisfy Assumptions \ref{asm24} and
\ref{asm27}.
For more choices of $w_t$, we refer to \citet{r24} and \citet{r30}.
We first state the strong convergence of $\hat{\theta}_{sn}$
in the following theorem and its proof is given in the \hyperref[app]{Appendix}.
%
%th2.1 #&#
\begin{theorem}\label{theorem21}
Suppose $\eta_t$ has a median zero with $E|\eta_t|=1$. If
Assumptions \ref{asm21}--\ref{asm25} hold, then
\[
\hat{\theta}_{sn}\rightarrow\theta_0 \qquad\mbox{a.s., as } n\to
\infty.
\]
\end{theorem}

To study the rate of convergence of $\hat{\theta}_{sn}$, we
reparameterize the weighted log-likelihood function (\ref{23}) as
follows:
\[
L_n(u)\equiv nL_{sn}(\theta_0+u)-nL_{sn}(\theta_0),
\]
where $u\in\Lambda\equiv\{u=(u_1',u_2')':u+\theta_0\in\Theta\}$. Let
$\hat{u}_{n}=\hat{\theta}_{sn}-\theta_0$. Then, $\hat{u}_{n}$ is the
minimizer of $L_n(u)$ on $\Lambda$. Furthermore, we have
%
%e2.5 ###
%e2.5 #&#
\begin{equation}\label{25}
L_n(u)=\sum_{t=1}^nw_tA_t(u)+\sum_{t=1}^nw_tB_t(u)+\sum_{t=1}^nw_tC_t(u),
\end{equation}
where
\begin{eqnarray*}
A_t(u)&=&\f{1}{\sqrt{h_t(\theta_0)}}[|\e_t(\gamma_0+u_1)|-
|\e_t(\gamma_0)|],\\
B_t(u)&=&\log\sqrt{h_t(\theta_0+u)}-\log\sqrt{h_t(\theta_0)}
+\f{|\e_t(\gamma_0)|}{\sqrt{h_t(\theta_0+u)}}
-\f{|\e_t(\gamma_0)|}{\sqrt{h_t(\theta_0)}},\\
C_t(u)&=&\biggl[\f{1}{\sqrt{h_t(\theta_0+u)}}
-\f{1}{\sqrt{h_t(\theta_0)}}\biggr][|\e_t(\gamma_0+u_1)|-|\e_t(\gamma_0)|].
\end{eqnarray*}

Let $I(\cdot)$ be the indicator function. Using the identity
%
%e2.6 ###
%e2.6 #&#
\begin{eqnarray} \label{26}
|x-y|-|x|&=&-y[I(x>0)-I(x<0)]\nonumber\\[-8pt]\\[-8pt]
&&{}+2\int_{0}^y [I(x\leq s)-I(x\leq0)]\,ds\nonumber
\end{eqnarray}
for $x\neq0$, we can show that
%
%e2.7 ###
%e2.7 #&#
\begin{equation} \label{27}
A_t(u)=q_t(u)[I(\eta_t>0)-I(\eta_t<0)]+2\int_{0}^{-q_t(u)}
X_t(s)\,ds,
\end{equation}
where $X_t(s)=I(\eta_t\leq s)-I(\eta_t\leq0)$,
$q_t(u)=q_{1t}(u)+q_{2t}(u)$ with
\[
q_{1t}(u)=\f{u'}{\sqrt{h_t(\theta_0)}}\,\f{\p\e_t(\gamma_0)}{\p\theta}
\quad\mbox{and}\quad
q_{2t}(u)=\f{u'}{2\sqrt{h_t(\theta_0)}}\,\f{\p^2\e_t(\xi^*)}{\p\theta\,\p
\theta'}u,
\]
and $\xi^*$ lies between $\gamma_0$ and $\gamma_0+u_1$. Moreover,
let $\mathcal{F}_{t}=\sigma\{\eta_k: k\leq t\}$ and
\[
\xi_t(u)=2w_t\int_{0}^{-q_{1t}(u)} X_t(s) \,ds.
\]
Then, from (\ref{27}), we have
%
%e2.8 ###
%e2.8 #&#
\begin{equation}\label{28}
\sum_{t=1}^n w_tA_t(u) = u' T_{1n} +
\Pi_{1n}(u)+\Pi_{2n}(u)+\Pi_{3n}(u),
\end{equation}
where
\begin{eqnarray*}
T_{1n}&=&\sum_{t=1}^{n}
\f{w_t}{\sqrt{h_t(\theta_0)}}\,\f{\p\e_t(\gamma_0)}{\p\theta}
[I(\eta_t>0)-I(\eta_t<0)],\\
\Pi_{1n}(u)&=&\sum_{t=1}^n
\{\xi_t(u)-E[\xi_t(u)|\mathcal{F}_{t-1}]\},\\
\Pi_{2n}(u)&=&\sum_{t=1}^nE[\xi_t(u)|\mathcal{F}_{t-1}],\\
\Pi_{3n}(u)&=&\sum_{t=1}^n w_tq_{2t}(u)[I(\eta_t>0)-I(\eta_t<0)]\\
&&{} +2 \sum_{t=1}^n w_t\int_{-q_{1t}(u)}^{-q_t(u)} X_t(s) \,ds.
\end{eqnarray*}

By Taylor's expansion, we can see that
%
%e2.9 ###
%e2.9 #&#
\begin{equation}\label{29}
\sum_{t=1}^n w_tB_t(u)=u'T_{2n}+\Pi_{4n}(u)+\Pi_{5n}(u),
\end{equation}
where
\begin{eqnarray*}
T_{2n}&=&\sum_{t=1}^n \f{w_t}{2h_t(\theta_0)}\,\f{\p
h_t(\theta_0)}{\p\theta}(1-|\eta_t|),\\
\Pi_{4n}(u)&=&u'\sum_{t=1}^n w_t
\biggl(\f{3}{8}\biggl|\f{\e_t(\gamma_0)}{\sqrt{h_t(\zeta^*)}}\biggr|
-\f{1}{4}\biggr)\f{1}{h_t^2(\zeta^*)} \,\f{\p h_t(\zeta^*)}{\p\theta}\,
\f{\p
h_t(\zeta^*)}{\p\theta'}u,\\
\Pi_{5n}(u)&=&u'\sum_{t=1}^n
w_t\biggl(\f{1}{4}-\f{1}{4}\biggl|\f{\e_t(\gamma_0)}{\sqrt{h_t(\zeta^*)}}\biggr|\biggr)
\f{1}{h_t(\zeta^*)}\,\f{\p^2 h_t(\zeta^*)}{\p\theta\,\p\theta'}u,
\end{eqnarray*}
and $\zeta^*$ lies between $\theta_0$ and $\theta_0+u$.

We further need one assumption and three lemmas. The first lemma is
directly from the central limit theorem for a martingale difference
sequence. The second- and third-lemmas give the expansions of
$\Pi_{in}(u)$ for $i=1,\ldots,5$ and $\sum_{t=1}^n C_t(u)$. The key
technical argument is for the second lemma for which we use the
bracketing method in \citet{r32}.
%
%as2.6 #&#
\begin{asm}\label{asm26}
$\eta_t$ has zero median with $E|\eta_t|=1$ and a
continuous density
function $g(x)$ satisfying $g(0)>0$ and $\sup_{x\in R}g(x)<\infty$.
\end{asm}
%
%le2.1 #&#
\begin{lem}\label{lemma21}
Let $T_n=T_{1n}+T_{2n}$. If Assumptions \ref{asm21}--\ref{asm26} hold, then
\[
\f{1}{\sqrt{n}}T_n\rightarrow_d N(0,\Omega_0) \qquad\mbox{as }
n\to\infty,
\]
where $\to_d$ denotes the convergence in distribution and
\[
\Omega_0=E\biggl(\f{w_t^{2}}{h_t(\theta_0)}\,\f{\p\e_t(\gamma_0)}{\p\theta}\,\f{\p
\e_t(\gamma_0)}{\p\theta'}\biggr)
+\f{E\eta_t^{2}-1}{4}E\biggl(\f{w_t^{2}}{h_t^{2}(\theta_0)}\,\f{\p
h_t(\theta_0)}{\p\theta}\,\f{\p h_t(\theta_0)}{\p\theta'}\biggr).
\]
\end{lem}
%
%le2.2 #&#
\begin{lem}\label{lemma22}
If Assumptions \ref{asm21}--\ref{asm26} hold, then for any sequence of random
variables $u_n$ such that $u_n=o_p(1)$, it follows that
\[
\Pi_{1n}(u_n)=o_p\bigl(\sqrt{n}\|u_n\|+n\|u_n\|^2\bigr),
\]
where $o_p(\cdot)\to0$ in probability as $n\to\infty$.
\end{lem}
%
%le2.3 #&#
\begin{lem}\label{lemma23}
If Assumptions \ref{asm21}--\ref{asm26} hold, then for any sequence of random
variables $u_n$ such that $u_n=o_p(1)$, it follows that:
\begin{eqnarray*}
\mbox{\textup{(i)}\quad\hspace*{10pt}} \Pi_{2n}(u_n)&=&\bigl(\sqrt{n}u_n\bigr)'\Sigma_1\bigl(\sqrt{n}u_n\bigr)+o_p(n\|u_n\|
^2),\\
\mbox{\textup{(ii)}\quad\hspace*{10pt}} \Pi_{3n}(u_n)&=&o_p(n\|u_n\|^2),\\
\mbox{\textup{(iii)}\quad\hspace*{10pt}} \Pi_{4n}(u_n)&=&\bigl(\sqrt{n}u_n\bigr)'\Sigma_2\bigl(\sqrt{n}u_n\bigr)+o_p(n\|u_n\|
^2),\\
\mbox{\textup{(iv)}\quad\hspace*{10pt}} \Pi_{5n}(u_n)&=&o_p(n\|u_n\|^2),\\
\mbox{\textup{(v)}\quad} \sum_{t=1}^n C_t(u_n)&=&o_p(n\|u_n\|^2),
\end{eqnarray*}
where
\[
\Sigma_1=g(0)E\biggl(\f{w_t}{h_t(\theta_0)}\,\f{\p\e_t(\gamma_0)}{\p\theta}\,\f{\p
\e_t(\gamma_0)}{\p\theta'}\biggr)
\]
and
\[
\Sigma_2=\frac{1}{8}E\biggl(\f{w_t}{h_t^2(\theta_0)}\,\f{\p
h_t(\theta_0)}{\p\theta}\,\f{\p h_t(\theta_0)}{\p\theta'}\biggr).
\]
\end{lem}

The proofs of Lemmas \ref{lemma22} and \ref{lemma23} are given in Section
\ref{sec6}. We now can state
one main result as follows:
%
%th2.2 #&#
\begin{theorem}\label{theorem22}
If Assumptions \ref{asm21}--\ref{asm26} hold, then:
\begin{eqnarray*}
\mbox{\textup{(i)}\quad} \sqrt{n}(\hat{\theta}_{sn}-\theta_0)&=&O_p(1),\\
\mbox{\textup{(ii)}\quad} \sqrt{n}(\hat{\theta}_{sn}-\theta_0)&\to_d& N\bigl(0,
\tfrac{1}{4}\Sigma_0^{-1}\Omega_0\Sigma_0^{-1}\bigr) \qquad\mbox{as }
n\to\infty,
\end{eqnarray*}
where $\Sigma_0=\Sigma_1+\Sigma_2$.
\end{theorem}

\begin{pf}
(i) First, we have $\hat{u}_n=o_p(1)$ by Theorem \ref{theorem21}. Furthermore,
by~(\ref{25}), (\ref{28}) and (\ref{29}) and Lemmas \ref{lemma22} and \ref
{lemma23}, we have
%
%e2.10 ###
%e2.10 #&#
\begin{equation} \label{210}
\qquad L_n(\hat{u}_n)=\hat{u}_n' T_n + \bigl(\sqrt{n}\hat{u}_n\bigr)'\Sigma_0
\bigl(\sqrt{n}\hat{u}_n\bigr)+ o_p\bigl(\sqrt{n}\|\hat{u}_n\|+n\|\hat{u}_n\|^2\bigr).
\end{equation}
Let $\lambda_{\min}>0$ be the minimum eigenvalue of $\Sigma_0$. Then
\[
L_n(\hat{u}_n)\geq-\bigl\|\sqrt{n}\hat{u}_n\bigr\|
\biggl[\biggl\|\f{1}{\sqrt{n}}T_n\biggr\|+o_p(1)\biggr]+n\|\hat{u}_n\|^2[\lambda_{\min}+o_p(1)].
\]
Note that $L_n(\hat{u}_{n})\leq0$. By the previous inequality, it
follows that
%
%e2.11 ###
%e2.11 #&#
\begin{equation}\label{211}
\sqrt{n}\| \hat{u}_n\|\le[\lambda_{\min}+o_p(1)]^{-1}
\biggl[\biggl\|\f{1}{\sqrt{n}}T_n\biggr\|+o_p(1)\biggr]=O_p(1),
\end{equation}
where the last step holds by Lemma \ref{lemma21}. Thus, (i) holds.

(ii) Let $u^*_n=-\Sigma_0^{-1}T_n/2n$. Then, by Lemma \ref{lemma21},
we have
\[
\sqrt{n}u^*_n\rightarrow_d
N\bigl(0,\tfrac{1}{4}\Sigma_0^{-1}\Omega_0\Sigma_0^{-1}\bigr)
\qquad\mbox{as } n\to\infty.
\]
Hence, it is sufficient to show that
$\sqrt{n}\hat{u}_n-\sqrt{n}u_n^*=o_p(1)$. By (\ref{210}) and~(\ref{211}), we have
\begin{eqnarray*}
L_n(\hat{u}_n)&=&\bigl(\sqrt{n}\hat{u}_n\bigr)'
\frac{1}{\sqrt{n}}T_n+\bigl(\sqrt{n}\hat
{u}_n\bigr)'\Sigma_0
\bigl(\sqrt{n}\hat{u}_n\bigr)+o_p(1)\\
&=&\bigl(\sqrt{n}\hat{u}_n\bigr)'\Sigma_0
\bigl(\sqrt{n}\hat{u}_n\bigr)-2\bigl(\sqrt{n}\hat{u}_n\bigr)'\Sigma_0
\bigl(\sqrt{n}u^*_n\bigr)+o_p(1).
\end{eqnarray*}
Note that (\ref{210}) still holds when $\hat{u}_n$ is replaced by
$u^*_n$. Thus,
\begin{eqnarray*}
L_n(u_n^*)&=&\bigl(\sqrt{n}u_n^*\bigr)'\frac{1}{\sqrt{n}}T_n+\bigl(\sqrt{n}u_n^*\bigr)'
\Sigma_0
\bigl(\sqrt{n}u_n^*\bigr)+o_p(1) \\
&=&-\bigl(\sqrt{n}u_n^*\bigr)'\Sigma_0 \bigl(\sqrt{n}u_n^*\bigr)+o_p(1).
\end{eqnarray*}
By the previous two equations, it follows that
%
%e2.12 ###
%e2.12 #&#
\begin{eqnarray}\label{212}
\qquad L_n(\hat{u}_n)-L_n(u^*_n)&=&\bigl(\sqrt{n}\hat{u}_n-\sqrt{n}u_n^*\bigr)'\Sigma_0
\bigl(\sqrt{n}\hat{u}_n-\sqrt{n}u_n^*\bigr)+o_p(1)
\nonumber\\[-8pt]\\[-8pt]
\qquad&\ge&\lambda_{\min}\bigl\|\sqrt{n}\hat{u}_n-\sqrt{n}u_n^*\bigr\|^{2}+o_{p}(1).\nonumber
\end{eqnarray}
Since
$L_n(\hat{u}_n)-L_n(u^*_n)=n[L_{sn}(\theta_{0}+\hat{u}_n)-L_{sn}(\theta
_{0}+u^*_n)]
\le0$ a.s., by (\ref{212}), we have
$\|\sqrt{n}\hat{u}_n-\sqrt{n}u_n^*\|=o_{p}(1)$. This completes the
proof.
\end{pf}
%
%re2.1 #&#
\begin{rem}
When $w_{t}\equiv1$, the limiting distribution in Theorem \ref
{theorem22} is the
same as that in \citet{r23}. When $r=s=0$ (ARMA model), it
reduces to the case in \citet{r30} and \citet{r36}.
In general, it is not easy to compare the asymptotic efficiency of
the self-weighted QMELE and the self-weight QMLE in \citet{r25}.
However, for the pure ARCH model, a~formal comparison of these two
estimators is given in Section \ref{sec3}. For the general ARMA--GARCH model,
a comparison based on simulation is given in Section \ref{sec4}.
\end{rem}

In practice, the initial values $Y_{0}$ are unknown, and have to be
replaced by some constants. Let $\tilde{\e}_{t}(\theta)$,
$\tilde{h}_{t}(\theta)$ and $\tilde{w}_{t}$ be $\e_{t}(\theta)$,
$h_{t}(\theta)$ and $w_{t}$, respectively, when $Y_{0}$ are
constants not depending on parameters. Usually, $Y_{0}$ are taken to
be zeros. The objective function (\ref{23}) is modified as
\[
\tilde{L}_{sn}(\theta)=\frac{1}{n}\sum_{t=1}^{n}
\tilde{w}_{t}\tilde{l}_{t}(\theta) \quad\mbox{and}\quad
\tilde{l}_{t}(\theta)=\log\sqrt{\tilde{h}_t(\theta)}+\f{|\tilde{\e
}_t(\gamma)|}{\sqrt{\tilde{h}_t(\theta)}}.
\]
To make the initial values $Y_{0}$ ignorable, we need the following
assumption.
%
%as2.7 #&#
\begin{asm} \label{asm27}
$E|w_{t}-\tilde{w}_{t}|^{\iota_{0}/4}=O(t^{-2})$, where
$\iota_{0}=\min\{\iota,1\}$.
\end{asm}

Let $\tilde{\theta}_{sn}$ be the minimizer of
$\tilde{L}_{sn}(\theta)$, that is,
\[
\tilde{\theta}_{sn}=\mathop{\arg\min}_{\Theta}\tilde{L}_{sn}(\theta).
\]
Theorem \ref{theorem23} below shows that $\tilde{\theta}_{sn}$ and
$\hat{\theta}_{sn}$ have the same limiting property. Its proof is
straightforward and can be found in \citet{r34}.
%
%th2.3 #&#
\begin{theorem}\label{theorem23}
Suppose that Assumption \ref{asm27} holds. Then, as $n\to\infty$,
\begin{eqnarray*}
&&\mbox{\hphantom{i}\textup{(i)}\quad}\mbox{if the assumptions of Theorem \ref{theorem21} hold}\\
&&\qquad\qquad\tilde
{\theta}_{sn}\to\theta_{0} \qquad\mbox{a.s.},\\
&&\mbox{\textup{(ii)}\quad}\mbox{if the assumptions of Theorem \ref{theorem22} hold}\\
&&\qquad\qquad\sqrt{n}(\tilde{\theta}_{sn}-\theta_0)\rightarrow_d
N\bigl(0,\tfrac{1}{4}\Sigma_0^{-1}\Omega_0\Sigma_0^{-1}\bigr).
\end{eqnarray*}
\end{theorem}

%s3 ###
%se3 #&#
\section{Local QMELE}\label{sec3}
The self-weighted QMELE in Section \ref{sec2} reduces the moment condition of
$\e_t$, but it may not be efficient. In this section, we propose a
local QMELE based on the self-weighted QMELE and derive its
asymptotic property. For some special cases, a formal comparison of
the local QMELE and the self-weighted QMELE is given.\vspace*{1pt}

Using $\hat{\theta}_{sn}$ in Theorem \ref{theorem22} as an initial
estimator of
$\theta_0$, we obtain the local QMELE $\hat{\theta}_n$ through the
following one-step iteration:
%
%e3.1 ###
%e3.1 #&#
\begin{equation}\label{31}
\hat{\theta}_n=\hat{\theta}_{sn}-[2\Sigma_n^{*}(\hat{\theta
}_{sn})]^{-1}T_n^{*}(\hat{\theta}_{sn}),
\end{equation}
where
\begin{eqnarray*}
\Sigma_n^{*}(\theta)&=&\sum_{t=1}^n\biggl\{\f{g(0)}{h_t(\theta)}\,\f{\p\e
_t(\gamma)}{\p\theta}\,\f{\p\e_t(\gamma)}{\p\theta'}
+ \f{1}{8h_t^2(\theta)}\,\f{\p
h_t(\theta)}{\p\theta}\,\f{\p h_t(\theta)}{\p\theta'}\biggr\},\\
T_n^{*}(\theta)&=&\sum_{t=1}^n \biggl\{\f{1}{\sqrt{h_t(\theta)}}\,\f{\p\e
_t(\gamma)}{\p\theta} \bigl[I\bigl(\eta_t(\theta)>0\bigr)-I\bigl(\eta_t(\theta)<0\bigr)\bigr] \\
&&\hspace*{87.3pt}{} +\f{1}{2h_t(\theta)}\,\f{\p
h_t(\theta)}{\p\theta}\bigl(1-|\eta_t(\theta)|\bigr)\biggr\}.
\end{eqnarray*}
In order to get the asymptotic normality of $\hat{\theta}_n$, we
need one more assumption as follows:
%
%as3.1 #&#
\begin{asm} \label{asm31}
$E\eta_t^2\sum_{i=1}^{r} \alpha_{0i}+\sum_{i=1}^{s} \beta_{0i}<1$ or
\[
E\eta_t^2\sum_{i=1}^{r} \alpha_{0i}+\sum_{i=1}^{s} \beta_{0i}=1
\]
with $\eta_t$ having a positive density on $R$ such that
$E|\eta_t|^{\tau}<\infty$ for all $\tau<\tau_0$ and
$E|\eta_t|^{\tau_0}=\infty$ for some $\tau_0\in(0,\infty]$.
\end{asm}

Under Assumption \ref{asm31}, there exists a unique strictly stationary
causal solution to GARCH model (\ref{12}); see \citet{r8}
and \citet{r1}. The condition $E\eta_t^2\sum_{i=1}^{r}
\alpha_{0i}+\sum_{i=1}^{s} \beta_{0i}<1$ is necessary and sufficient
for $E\e_t^2<\infty$ under which model (\ref{12}) has a finite variance.
When $E\eta_t^2\sum_{i=1}^{r} \alpha_{0i}+\sum_{i=1}^{s}
\beta_{0i}=1$, model (\ref{12}) is called IGARCH model. The IGARCH model
has an infinite variance, but $E|\e_t|^{2\iota}<\infty$ for all
$\iota\in(0,1)$ under Assumption \ref{asm31}; see \citet{r25}. Assumption
\ref{asm31} is crucial for the ARMA--IGARCH model. From Figure \ref
{figure2} in Section
\ref{sec2}, we can see that the parameter\vspace*{1pt} region specified in Assumption \ref{asm31}
is much bigger than that for $E|\e_t|^{3}<\infty$ which is required
for the asymptotic normality of the global QMELE. Now, we give one
lemma as follows and its proof is straightforward and can be found
in \citet{r34}.
%
%le3.1 #&#
\begin{lem}\label{lemma31}
If Assumptions \ref{asm21}--\ref{asm23}, \ref{asm26} and \ref{asm31}
hold, then for any sequence of
random variables $\theta_n$ such that
$\sqrt{n}(\theta_n-\theta_0)=O_p(1)$, it follows that:
\begin{eqnarray*}
\mbox{\textup{(i)}\quad} \f{1}{n}[T_n^{*}(\theta_{n})-T_n^{*}(\theta_0)]&=&[2\Sigma
+o_p(1)](\theta_n-\theta_0)+o_p\biggl(\f{1}{\sqrt{n}}\biggr),\\
\mbox{\textup{(ii)}\hspace*{50pt}\quad} \f{1}{n}\Sigma_n^{*}(\theta_{n})&=&\Sigma+o_p(1),\\
\mbox{\textup{(iii)}\hspace*{42pt}\quad} \f{1}{\sqrt{n}}T_n^{*}(\theta_0)&\to_d&
N(0,\Omega) \qquad\mbox{as } n\to\infty,
\end{eqnarray*}
where
\begin{eqnarray*}
\Omega&=&E\biggl(\f{1}{h_t(\theta_0)}\,\f{\p\e_t(\gamma_0)}{\p\theta}\,\f{\p\e
_t(\gamma_0)}{\p\theta'}\biggr)
+\f{E\eta_t^{2}-1}{4}E\biggl(\f{1}{h_t^{2}(\theta_0)}\,\f{\p
h_t(\theta_0)}{\p\theta}\,\f{\p h_t(\theta_0)}{\p\theta'}\biggr),\\
\Sigma&=&g(0)E\biggl(\f{1}{h_t(\theta_0)}\,\f{\p\e_t(\gamma_0)}{\p\theta}\,\f{\p\e
_t(\gamma_0)}{\p\theta'}\biggr)
+\f{1}{8}E\biggl(\f{1}{h_t^2(\theta_0)}\,\f{\p
h_t(\theta_0)}{\p\theta}\,\f{\p h_t(\theta_0)}{\p\theta'}\biggr).
\end{eqnarray*}
\end{lem}
%
%th3.1 #&#
\begin{theorem}\label{theorem31}
If the conditions in Lemma \ref{lemma31} are satisfied, then
\[
\sqrt{n}(\hat{\theta}_n-\theta_0)\rightarrow_d
N\bigl(0,\tfrac{1}{4}\Sigma^{-1}\Omega\Sigma^{-1}\bigr) \qquad\mbox{as }
n\to\infty.
\]
\end{theorem}

\begin{pf}
Note that $\sqrt{n}(\hat{\theta}_{sn}-\theta_{0})=O_p(1)$. By (\ref{31})
and Lemma \ref{lemma31}, we have that
\begin{eqnarray*}
\hat{\theta}_{n}&=&\hat{\theta}_{sn}-\biggl[\f{2}{n}\Sigma_{n}^{*}(\hat{\theta
}_{sn})\biggr]^{-1}
\biggl[\f{1}{n}T_n^{*}(\hat{\theta}_{sn})\biggr]\\
&=&\hat{\theta}_{sn}-[2\Sigma+o_p(1)]^{-1} \biggl\{\f{1}{n}T_n^{*}(\theta
_0)+[2\Sigma+o_p(1)](\hat{\theta}_{sn}-\theta_0)+o_p\biggl(\f{1}{\sqrt{n}}\biggr)\biggr\}
\\
&=&\theta_0+\f{\Sigma^{-1}T_n^{*}(\theta_0)}{2n}+o_p\biggl(\f{1}{\sqrt{n}}\biggr).
\end{eqnarray*}
It follows that
\[
\sqrt{n}(\hat{\theta}_n-\theta_0)=\f{\Sigma^{-1}T_n^{*}(\theta
_0)}{2\sqrt{n}}+o_p(1).
\]
By Lemma \ref{lemma31}(iii), we can see that the conclusion holds. This
completes the proof.
\end{pf}
%
%re3.1 #&#
\begin{rem}
In practice, by using $\tilde{\theta}_{sn}$ in Theorem \ref{theorem23}
as an
initial estimator of $\theta_{0}$, the local QMELE has to be
modified as follows:
\[
\hat{\theta}_n=\tilde{\theta}_{sn}-[2\tilde{\Sigma}_n^{*}(\tilde{\theta
}_{sn})]^{-1}\tilde{T}_n^{*}(\tilde{\theta}_{sn}),
\]
where $\tilde{\Sigma}_{n}^{*}(\theta)$ and
$\tilde{T}_{n}^{*}(\theta)$ are defined in the same way as
$\Sigma_{n}^{*}(\theta)$ and $T_{n}^{*}(\theta)$, respectively, with
$\e_{t}(\theta)$ and
$h_{t}(\theta)$ being replaced by $\tilde{\e}_{t}(\theta)$ and
$\tilde{h}_{t}(\theta)$. However, this does not affect the
asymptotic property of $\hat{\theta}_{n}$; see Theorem 4.3.2 in
\citet{r34}.
\end{rem}

We now compare the asymptotic efficiency of the local QMELE and the
self-weighted QMELE. First, we consider the pure ARMA model, that is,
model (\ref{11})--(\ref{12}) with $h_t$ being a constant. In this case,
\begin{eqnarray*}
\Omega_0&=&E(w_t^{2}X_{1t}X_{1t}'), \qquad
\Sigma_0=g(0)E(w_t X_{1t}X_{1t}'),\\
\Omega&=&E(X_{1t}X_{1t}')
\quad\mbox{and}\quad \Sigma=g(0)\Omega,
\end{eqnarray*}
where $X_{1t}=h_t^{-1/2}\p\e_t(\gamma_0)/\p\theta$. Let b and c be
two any $m$-dimensional constant vectors. Then,
\begin{eqnarray*}
c'\Sigma_0bb'\Sigma_0c&=&\bigl\{E\bigl[\bigl(c'\sqrt{g(0)}w_tX_{1t}\bigr)\bigl(\sqrt
{g(0)}X_{1t}'b\bigr)\bigr]\bigr\}^2\\
&\leq& E\bigl(c'\sqrt{g(0)}w_tX_{1t}\bigr)^2E\bigl(\sqrt{g(0)}X_{1t}'b\bigr)^2\\
&=&[c'g(0)\Omega_0c][b'\Sigma b]=c'[g(0)\Omega_0b'\Sigma b]c.
\end{eqnarray*}
Thus, $g(0)\Omega_0b'\Sigma b'-\Sigma_0bb'\Sigma_0\geq0$ (a
positive semi-definite matrix) and hence $b'\Sigma_0\Omega_0^{-1}\Sigma
_0b=\operatorname{tr}(\Omega_0^{-1/2}\Sigma_0bb'\Sigma_0\Omega_0^{-1/2})\leq
\operatorname{tr}(g(0)b'\Sigma b)=g(0)b'\Sigma b$. It follows that
$\Sigma_0^{-1}\Omega_0\Sigma_0^{-1}\geq
[g(0)\Sigma]^{-1}=\Sigma^{-1}\Omega\Sigma^{-1}$. Thus, the local
QMELE is more efficient than the self-weighted QMELE. Similarly, we
can show that the local QMELE is more efficient than the
self-weighted QMELE for the pure GARCH model.

For the general model (\ref{11})--(\ref{12}), it is not easy to compare the
asymptotic efficiency of the self-weighted QMELE and the local
QMELE. However, when $\eta_t\sim \operatorname{Laplace}(0,1)$, we have
\begin{eqnarray*}
\Sigma_0&=&E\biggl(\f{w_t}{2}X_{1t}X_{1t}'+\f{w_t}{8}X_{2t}X_{2t}'\biggr),\\
\Omega_0&=&E\biggl(w_t^2X_{1t}X_{1t}'+\f{w_t^2}{4}X_{2t}X_{2t}'\biggr),\\
\Sigma&=&E\bigl(\tfrac{1}{2}X_{1t}X_{1t}'+\tfrac{1}{8}X_{2t}X_{2t}'\bigr)
\quad\mbox{and}\quad \Omega=2\Sigma,
\end{eqnarray*}
where $X_{2t}=h_t^{-1}\p h_t(\theta_0)/\p\theta$. Then, it is easy
to see that
\begin{eqnarray*}
&&c'\Sigma_0bb'\Sigma_0c\\
&&\qquad=\{
E[(c'2^{-1/4}w_tX_{1t})(2^{-3/4}X_{1t}'b)+(c'2^{-5/4}w_tX_{2t})(2^{-7/4}X_{2t}'b)]\}
^2\\
&&\qquad\leq\bigl\{\sqrt{E(c'2^{-1/4}w_tX_{1t})^2E(2^{-3/4}X_{1t}'b)^2}\\
&&\qquad\quad\hspace*{4pt}{}+\sqrt
{E(c'2^{-5/4}w_tX_{2t})^2E(2^{-7/4}X_{2t}'b)^2}\bigr\}^2\\
&&\qquad\leq
[E(c'2^{-1/4}w_tX_{1t})^2+E(c'2^{-5/4}w_tX_{2t})^2]\\
&&\qquad\quad{}\times[E(2^{-3/4}X_{1t}'b)^2+E(2^{-7/4}X_{2t}'b)^2]\\
&&\qquad=[c'2^{-1/2}\Omega_0c][b'2^{-1/2}\Sigma
b]=c'[2^{-1}\Omega_0b'\Sigma b]c.
\end{eqnarray*}
Thus, $2^{-1}\Omega_0b'\Sigma b'-\Sigma_0bb'\Sigma_0\geq0$ and hence
$b'\Sigma_0\Omega_0^{-1}\Sigma_0b=\operatorname{tr}(\Omega_0^{-1/2}
\Sigma_0bb'\times\Sigma
_0\Omega_0^{-1/2})
\leq \operatorname{tr}(2^{-1}b'\Sigma b)=2^{-1}b'\Sigma b$.
It follows that $\Sigma_0^{-1}\Omega_0\Sigma_0^{-1}\geq
2\Sigma^{-1}=\break\Sigma^{-1}\Omega\Sigma^{-1}$. Thus, the local QMELE is
more efficient than the global self-weighted QMELE.

In the end, we compare the asymptotic efficiency of the
self-weighted QMELE and the self-weighted QMLE in \citet{r25} for
the pure ARCH model, when $E\eta_t^4<\infty$. We reparametrize model
(\ref{12}) when $s=0$ as follows:
%
%e3.2 ###
%e3.2 #&#
\begin{equation}\label{32}
y_t=\eta_t^{*}\sqrt{h_t^{*}} \quad\mbox{and}\quad
h_t^{*}=\alpha_0^{*}+\sum_{i=1}^{r} \alpha_i^{*}y_{t-i}^2,
\end{equation}
where $\eta_t^{*}=\eta_t/\sqrt{E\eta_t^2}$, $h_t^{*}=(E\eta_t^2)h_t$
and
$\theta^{*}=(\alpha_0^{*},\alpha_1^{*},\ldots,\alpha_r^{*})'=(E\eta
_t^2)\theta$.\break Let~$\tilde{\theta}_{sn}^{*}$ be the self-weighted QMLE of the true
parameter, $\theta_0^{*}$, in model~(\ref{32}). Then,
$\tilde{\theta}_{sn}=\tilde{\theta}_{sn}^{*}/E\eta_t^2$ is the
self-weighted QMLE of $\theta_0$, and its asymptotic covariance is
\[
\Gamma_1=\kappa_1
[E(w_tX_{2t}X_{2t}')]^{-1}E(w_t^2X_{2t}X_{2t}')[E(w_tX_{2t}X_{2t}')]^{-1},
\]
where $\kappa_1=E\eta_t^4/(E\eta_t^2)^2-1$. By Theorem \ref{theorem22}, the
asymptotic variance of the self-weighted QMELE is
\[
\Gamma_2=\kappa_2
[E(w_tX_{2t}X_{2t}')]^{-1}E(w_t^2X_{2t}X_{2t}')[E(w_tX_{2t}X_{2t}')]^{-1},
\]
where $\kappa_2=4(E\eta_t^2-1)$. When $\eta_t\sim \operatorname{Laplace}(0,1)$,
$\kappa_1=5$ and $\kappa_2=4$. Thus, $\Gamma_1>\Gamma_2$, meaning
that the self-weighted QMELE is more efficient than the
self-weighted QMLE. When
$\eta_{t}=\tilde{\eta}_t/E|\tilde{\eta}_t|$, with $\tilde{\eta}_t$
having the following mixing normal density:
\[
f(x)=(1-\varepsilon)\phi(x)+\f{\varepsilon}{\tau}\phi\biggl(\f{x}{\tau}\biggr),
\]
we have $E|\eta_{t}|=1$,
\[
E\eta^2_t=\f{\pi(1-\e+\e\tau^2)}{2(1-\e+\e\tau)^2}
\]
and
\[
E\eta_t^4=\f{3\pi(1-\e+\e\tau^4)}{2(1-\e+\e\tau)^{2}(1-\e+\e\tau^2)},
\]
where $\phi(x)$ is the pdf of standard normal, $0\leq\e\leq1$ and
$\tau>0$. The asymptotic efficiencies of the self-weighted QMELE and
the self-weighted QMLE depend on $\e$ and $\tau$. For example, when
$\e=1$ and $\tau=\sqrt{\pi/2}$, we have $\kappa_1=(6-\pi)/\pi$ and
$\kappa_2=2\pi-4$, and hence the self-weighted QMLE is more
efficient than the self-weighted QMELE since $\Gamma_1<\Gamma_2$.
When $\e=0.99$ and $\tau=0.1$, we have $\kappa_1=28.1$ and
$\kappa_2=6.5$, and hence the self-weighted QMELE is more efficient
than the self-weighted QMLE since $\Gamma_1>\Gamma_2$.

%s4 ###
%se4 #&#
\section{Simulation}\label{sec4}
In this section, we compare the performance of the global
self-weighted QMELE ($\hat{\theta}_{sn}$), the global self-weighted
QMLE ($\bar{\theta}_{sn}$), the local QMELE $(\hat{\theta}_n)$ and
the local QMLE $(\bar{\theta}_n)$. The following AR(1)--GARCH$(1,1)$
model is used to generate data samples:
%
%e4.1 ###
%e4.1 #&#
\begin{eqnarray}\label{41}
y_t&=&\mu+\phi_1 y_{t-1}+\e_t, \nonumber\\[-8pt]\\[-8pt]
\e_t&=&\eta_t\sqrt{h_t} \quad\mbox{and}\quad
h_t=\alpha_0+\alpha_1\e_{t-1}^2+\beta_1h_{t-1}.\nonumber
\end{eqnarray}
We set the sample size $n=1\mbox{,}000$ and use $1\mbox{,}000$ replications, and
study the cases when $\eta_t$ has $\operatorname{Laplace}(0,1)$, $N(0,1)$ and $t_3$
distribution. For the case with $E\e_{t}^{2}<\infty$ (i.e.,
$E\eta_{t}^{2} \alpha_{01}+\beta_{01}<1$), we take\vspace*{1pt}
$\theta_0=(0.0,0.5,0.1,0.18,0.4)$. For the IGARCH case (i.e.,
$E\eta_{t}^{2} \alpha_{01}+\beta_{01}=1$), we take
$\theta_0=(0.0,0.5,0.1,0.3,0.4)$ when $\eta_t\sim \operatorname{Laplace}(0,1)$,
$\theta_0=(0.0,0.5,0.1,0.6,0.4)$ when $\eta_t\sim N(0,1)$ and
$\theta_0=(0.0,0.5,0.1,0.2,0.4)$ when $\eta_t\sim t_3$. We
standardize the distribution of~$\eta_t$ to ensure that
$E|\eta_t|=1$ for the QMELE. Tables \ref{table1}--\ref{table3} list the
sample biases, the
sample standard deviations (SD)\vspace*{1pt} and the asymptotic standard
deviations (AD) of $\hat{\theta}_{sn}, \bar{\theta}_{sn}$,
$\hat{\theta}_n$ and $\bar{\theta}_{n}$. We choose $w_{t}$ as in
(\ref{24}) with $C$ being 90\% quantile of $\{y_{1},\ldots,y_{n}\}$ and
$y_{i}\equiv0$ for $i\leq0$. The ADs\vspace*{1pt} in Theorems~\ref{theorem22}
and \ref{theorem31} are estimated by
$\hat{\chi}_{sn}=1/4\hat{\Sigma}_{sn}^{-1}\hat{\Omega}_{sn}\hat{\Sigma
}_{sn}^{-1}$
and
$\hat{\chi}_{n}=1/4\hat{\Sigma}_{n}^{-1}\hat{\Omega}_{n}\hat{\Sigma}_{n}^{-1}$,
respectively, where
{\small\begin{eqnarray*}
\hat{\Sigma}_{sn}&=&\f{1}{n}\sum_{t=1}^{n}
\biggl\{\f{g(0)w_t}{h_t(\hat{\theta}_{sn})}\,\f{\p\e_t(\hat{\gamma}_{sn})}{\p
\theta}\,\f{\p\e_t(\hat{\gamma}_{sn})}{\p\theta'}
+\f{w_t}{8h_t^{2}(\hat{\theta}_{sn})}\,\f{\p
h_t(\hat{\theta}_{sn})}{\p\theta}\,\f{\p
h_t(\hat{\theta}_{sn})}{\p\theta'}\biggr\},\\
\hat{\Omega}_{sn}&=&\f{1}{n}\sum_{t=1}^{n}
\biggl\{\f{w_t^{2}}{h_t(\hat{\theta}_{sn})}\,\f{\p\e_t(\hat{\gamma}_{sn})}{\p
\theta}\,\f{\p\e_t(\hat{\gamma}_{sn})}{\p\theta'}
+\f{E\eta_t^{2}-1}{4}\f{w_t^{2}}{h_t^{2}(\hat{\theta}_{sn})}\,\f{\p
h_t(\hat{\theta}_{sn})}{\p\theta}\,\f{\p
h_t(\hat{\theta}_{sn})}{\p\theta'}\biggr\},\\
\hat{\Sigma}_{n}&=&\f{1}{n}\sum_{t=1}^{n}
\biggl\{\f{g(0)}{h_t(\hat{\theta}_{n})}\,\f{\p\e_t(\hat{\gamma}_{n})}{\p\theta
}\,\f{\p\e_t(\hat{\gamma}_{n})}{\p\theta'}
+\f{1}{8h_t^{2}(\hat{\theta}_{n})}\,\f{\p
h_t(\hat{\theta}_{n})}{\p\theta}\,\f{\p
h_t(\hat{\theta}_{n})}{\p\theta'}\biggr\},\\
\hat{\Omega}_{n}&=&\f{1}{n}\sum_{t=1}^{n}
\biggl\{\f{1}{h_t(\hat{\theta}_{n})}\,\f{\p\e_t(\hat{\gamma}_{n})}{\p\theta}\,\f
{\p\e_t(\hat{\gamma}_{n})}{\p\theta'}
+\f{E\eta_t^{2}-1}{4}\f{1}{h_t^{2}(\hat{\theta}_{n})}\,\f{\p
h_t(\hat{\theta}_{n})}{\p\theta}\,\f{\p
h_t(\hat{\theta}_{n})}{\p\theta'}\biggr\}.
\end{eqnarray*}}

%
%t1 ###
%ta1 #&#
\begin{table}
\tabcolsep=0pt
\caption{Estimators for model (\protect\ref{41}) when $\eta_t\sim
\operatorname{Laplace}(0,1)$} \label{table1}
{\fontsize{8.5pt}{11pt}\selectfont{\begin{tabular*}{\tablewidth}
{@{\extracolsep{\fill}}ld{2.4}d{2.4}d{1.4}d{1.4}d{2.4}d{1.4}d{2.4}d{1.4}d{2.4}d{2.4}@{}}
\hline
&\multicolumn{5}{c}{$\bolds{\theta_0=(0.0, 0.5, 0.1, 0.18, 0.4)}$}
&\multicolumn{5}{c@{}}{$\bolds{\theta_0=(0.0, 0.5, 0.1, 0.3, 0.4)}$}\\[-4pt]
& \multicolumn{5}{c}{\hrulefill} & \multicolumn{5}{c}{\hrulefill}\\
&\multicolumn{5}{c}{\textbf{Self-weighted QMELE (}$\bolds{\hat{\theta
}_{sn}}$\textbf{)}} & \multicolumn{5}{c@{}}{\textbf{Self-weighted QMELE (}$\bolds{\hat{\theta
}_{sn}}$\textbf{)}}\\[-4pt]
& \multicolumn{5}{c}{\hrulefill} & \multicolumn{5}{c}{\hrulefill}\\
&\multicolumn{1}{c}{$\bolds{\hat{\mu}_{sn}}$}&\multicolumn{1}{c}{$\bolds{\hat{\phi}_{1sn}}$}&\multicolumn{1}{c}{$\bolds{\hat{\alpha}_{0sn}}$}
&\multicolumn{1}{c}{$\bolds{\hat{\alpha}_{1sn}}$}&\multicolumn{1}{c}{$\bolds{\hat{\beta}_{1sn}}$}&\multicolumn{1}{c}{$\bolds{\hat{\mu}_{sn}}$}
&\multicolumn{1}{c}{$\bolds{\hat{\phi}_{1sn}}$}&\multicolumn{1}{c}{$\bolds{\hat
{\alpha}_{0sn}}$} & \multicolumn{1}{c}{$\bolds{\hat{\alpha}_{1sn}}$}&\multicolumn{1}{c@{}}{$\bolds{\hat{\beta}_{1sn}}$}\\
\hline
Bias&0.0004&-0.0023&0.0034&0.0078&-0.0154&0.0003&-0.0049&0.0031&0.0054&-0.0068\\
SD&0.0172&0.0317&0.0274&0.0548&0.1125&0.0195&0.0318&0.0219&0.0640&0.0673\\
AD&0.0166&0.0304&0.0255&0.0540&0.1061&0.0192&0.0311&0.0218&0.0624&0.0664\\
\hline
%&\multicolumn{5}{c|}{}&\multicolumn{5}{c}{}\\
&\multicolumn{5}{c}{\textbf{Local QMELE (}$\bolds{\hat{\theta}_n}$\textbf{)}}
&\multicolumn{5}{c}{\textbf{Local QMELE (}$\bolds{\hat{\theta}_n}$\textbf{)}}\\[-4pt]
& \multicolumn{5}{c}{\hrulefill} & \multicolumn{5}{c}{\hrulefill}\\
&\multicolumn{1}{c}{$\bolds{\hat{\mu}_n}$}&\multicolumn{1}{c}{$\bolds{\hat{\phi}_{1n}}$}&\multicolumn{1}{c}{$\bolds{\hat{\alpha}_{0n}}$}
&\multicolumn{1}{c}{$\bolds{\hat{\alpha}_{1n}}$}&\multicolumn{1}{c}{$\bolds{\hat{\beta}_{1n}}$}&\multicolumn{1}{c}{$\bolds{\hat{\mu}_{n}}$}
&\multicolumn{1}{c}{$\bolds{\hat{\phi}_{1n}}$}&\multicolumn{1}{c}{$\bolds{\hat{\alpha}_{0n}}$}&\multicolumn{1}{c}{$\bolds{\hat{\alpha}_{1n}}$}
&\multicolumn{1}{c@{}}{$\bolds{\hat{\beta}_{1n}}$}\\
\hline
Bias&0.0008&-0.0019&0.0027&0.0002&-0.0094&0.0010&-0.0044&0.0024&-0.0008&-0.0025\\
SD&0.0170&0.0253&0.0249&0.0400&0.0989&0.0192&0.0261&0.0203&0.0502&0.0591\\
AD&0.0162&0.0245&0.0234&0.0407&0.0920&0.0190&0.0258&0.0206&0.0499&0.0591\\
\hline
&\multicolumn{5}{c}{\textbf{Self-weighted QMLE (}$\bolds{\bar{\theta
}_{sn}}$\textbf{)}}&\multicolumn{5}{c}{\textbf{Self-weighted QMLE
(}$\bolds{\bar{\theta}_{sn}}$\textbf{)}}\\[-4pt]
& \multicolumn{5}{c}{\hrulefill} & \multicolumn{5}{c}{\hrulefill}\\
&\multicolumn{1}{c}{$\bolds{\bar{\mu}_{sn}}$}&\multicolumn{1}{c}{$\bolds{\bar{\phi}_{1sn}}$}&\multicolumn{1}{c}{$\bolds{\bar{\alpha}_{0sn}}$}
&\multicolumn{1}{c}{$\bolds{\bar{\alpha}_{1sn}}$}&\multicolumn{1}{c}{$\bolds{\bar{\beta}_{1sn}}$}&\multicolumn{1}{c}{$\bolds{\bar{\mu}_{sn}}$}
&\multicolumn{1}{c}{$\bolds{\bar{\phi}_{1sn}}$}&\multicolumn{1}{c}{$\bolds{\bar{\alpha}_{0sn}}$}&\multicolumn{1}{c}{$\bolds{\bar{\alpha}_{1sn}}$}
&\multicolumn{1}{c@{}}{$\bolds{\bar{\beta}_{1sn}}$}\\
\hline
Bias&-0.0003&-0.0016&0.0041&0.0114&-0.0227&0.0005&-0.0039&0.0031&0.0104&-0.0127\\
SD&0.0243&0.0451&0.0301&0.0624&0.1237&0.0283&0.0458&0.0242&0.0750&0.0755\\
AD&0.0240&0.0443&0.0285&0.0607&0.1184&0.0283&0.0461&0.0243&0.0704&0.0741\\
\hline
&\multicolumn{5}{c}{\textbf{Local QMLE (}$\bolds{\bar{\theta}_n}$\textbf{)}}
&\multicolumn{5}{c}{\textbf{Local QMLE (}$\bolds{\bar{\theta}_n}$\textbf{)}}\\[-4pt]
& \multicolumn{5}{c}{\hrulefill} & \multicolumn{5}{c}{\hrulefill}\\
&\multicolumn{1}{c}{$\bolds{\bar{\mu}_n}$}&\multicolumn{1}{c}{$\bolds{\bar{\phi}_{1n}}$}&\multicolumn{1}{c}{$\bolds{\bar{\alpha}_{0n}}$}
&\multicolumn{1}{c}{$\bolds{\bar{\alpha}_{1n}}$}&\multicolumn{1}{c}{$\bolds{\bar{\beta}_{1n}}$}&\multicolumn{1}{c}{$\bolds{\bar{\mu}_{n}}$}
&\multicolumn{1}{c}{$\bolds{\bar{\phi}_{1n}}$}&\multicolumn{1}{c}{$\bolds{\bar{\alpha}_{0n}}$}&\multicolumn{1}{c}{$\bolds{\bar{\alpha}_{1n}}$}
&\multicolumn{1}{c@{}}{$\bolds{\bar{\beta}_{1n}}$}\\
\hline
Bias&0.0007&-0.0034&0.0026&0.0037&-0.0144&0.0022&-0.0045&0.0020&0.0044&-0.0081\\
SD&0.0243&0.0368&0.0279&0.0461&0.1115&0.0282&0.0377&0.0227&0.0579&0.0674\\
AD&0.0236&0.0361&0.0261&0.0459&0.1026&0.0281&0.0384&0.0230&0.0564&0.0659\\
\hline
\end{tabular*}}}
\vspace*{-3pt}
\end{table}

%
%
%t2 ###
%ta2 #&#
\begin{table}
\tabcolsep=0pt
\caption{Estimators for model (\protect\ref{41}) when $\eta_t\sim N(0,1)$}\label{table2}
\vspace*{-3pt}
{\fontsize{8.5pt}{11pt}\selectfont{\begin{tabular*}{\tablewidth}
{@{\extracolsep{\fill}}ld{2.4}d{2.4}d{1.4}d{1.4}d{2.4}d{2.4}d{2.4}d{1.4}d{2.4}d{2.4}@{}}
\hline
&\multicolumn{5}{c}{$\bolds{\theta_0=(0.0, 0.5, 0.1, 0.18, 0.4)}$}
&\multicolumn{5}{c@{}}{$\bolds{\theta_0=(0.0, 0.5, 0.1, 0.6, 0.4)}$}\\[-4pt]
& \multicolumn{5}{c}{\hrulefill} & \multicolumn{5}{c}{\hrulefill}\\
&\multicolumn{5}{c}{\textbf{Self-weighted QMELE (}$\bolds{\hat{\theta
}_{sn}}$\textbf{)}} & \multicolumn{5}{c@{}}{\textbf{Self-weighted QMELE (}$\bolds{\hat{\theta
}_{sn}}$\textbf{)}}\\[-4pt]
& \multicolumn{5}{c}{\hrulefill} & \multicolumn{5}{c}{\hrulefill}\\
&\multicolumn{1}{c}{$\bolds{\hat{\mu}_{sn}}$}&\multicolumn{1}{c}{$\bolds{\hat{\phi}_{1sn}}$}&\multicolumn{1}{c}{$\bolds{\hat{\alpha}_{0sn}}$}
&\multicolumn{1}{c}{$\bolds{\hat{\alpha}_{1sn}}$}&\multicolumn{1}{c}{$\bolds{\hat{\beta}_{1sn}}$}&\multicolumn{1}{c}{$\bolds{\hat{\mu}_{sn}}$}
&\multicolumn{1}{c}{$\bolds{\hat{\phi}_{1sn}}$}&\multicolumn{1}{c}{$\bolds{\hat
{\alpha}_{0sn}}$} & \multicolumn{1}{c}{$\bolds{\hat{\alpha}_{1sn}}$}&\multicolumn{1}{c@{}}{$\bolds{\hat{\beta}_{1sn}}$}\\
\hline
Bias&0.0003&-0.0042&0.0075&0.0065&-0.0372&-0.0008&-0.0034&0.0029&-0.0019&-0.0028\\
SD&0.0192&0.0457&0.0366&0.0600&0.1738&0.0255&0.0437&0.0204&0.0815&0.0512\\
AD&0.0189&0.0443&0.0379&0.0604&0.1812&0.0257&0.0424&0.0202&0.0809&0.0491\\
\hline
%&\multicolumn{5}{c|}{}&\multicolumn{5}{c}{}\\
&\multicolumn{5}{c}{\textbf{Local QMELE (}$\bolds{\hat{\theta}_n}$\textbf{)}}
&\multicolumn{5}{c}{\textbf{Local QMELE (}$\bolds{\hat{\theta}_n}$\textbf{)}}\\[-4pt]
& \multicolumn{5}{c}{\hrulefill} & \multicolumn{5}{c}{\hrulefill}\\
&\multicolumn{1}{c}{$\bolds{\hat{\mu}_n}$}&\multicolumn{1}{c}{$\bolds{\hat{\phi}_{1n}}$}&\multicolumn{1}{c}{$\bolds{\hat{\alpha}_{0n}}$}
&\multicolumn{1}{c}{$\bolds{\hat{\alpha}_{1n}}$}&\multicolumn{1}{c}{$\bolds{\hat{\beta}_{1n}}$}&\multicolumn{1}{c}{$\bolds{\hat{\mu}_{n}}$}
&\multicolumn{1}{c}{$\bolds{\hat{\phi}_{1n}}$}&\multicolumn{1}{c}{$\bolds{\hat{\alpha}_{0n}}$}&\multicolumn{1}{c}{$\bolds{\hat{\alpha}_{1n}}$}
&\multicolumn{1}{c@{}}{$\bolds{\hat{\beta}_{1n}}$}\\
\hline
Bias&0.0006&-0.0051&0.0061&0.0019&-0.0268&0.0000&-0.0040&0.0029&-0.0048&-0.0015\\
SD&0.0184&0.0372&0.0357&0.0487&0.1674&0.0252&0.0364&0.0197&0.0671&0.0472\\
AD&0.0183&0.0370&0.0350&0.0488&0.1652&0.0252&0.0359&0.0194&0.0685&0.0453\\
\hline
&\multicolumn{5}{c}{\textbf{Self-weighted QMLE (}$\bolds{\bar{\theta
}_{sn}}$\textbf{)}}&\multicolumn{5}{c}{\textbf{Self-weighted QMLE
(}$\bolds{\bar{\theta}_{sn}}$\textbf{)}}\\[-4pt]
& \multicolumn{5}{c}{\hrulefill} & \multicolumn{5}{c}{\hrulefill}\\
&\multicolumn{1}{c}{$\bolds{\bar{\mu}_{sn}}$}&\multicolumn{1}{c}{$\bolds{\bar{\phi}_{1sn}}$}&\multicolumn{1}{c}{$\bolds{\bar{\alpha}_{0sn}}$}
&\multicolumn{1}{c}{$\bolds{\bar{\alpha}_{1sn}}$}&\multicolumn{1}{c}{$\bolds{\bar{\beta}_{1sn}}$}&\multicolumn{1}{c}{$\bolds{\bar{\mu}_{sn}}$}
&\multicolumn{1}{c}{$\bolds{\bar{\phi}_{1sn}}$}&\multicolumn{1}{c}{$\bolds{\bar{\alpha}_{0sn}}$}&\multicolumn{1}{c}{$\bolds{\bar{\alpha}_{1sn}}$}
&\multicolumn{1}{c@{}}{$\bolds{\bar{\beta}_{1sn}}$}\\
\hline
Bias&-0.0001&-0.0039&0.0069&0.0089&-0.0361&-0.0006&-0.0016&0.0024&0.0027&-0.0045\\
SD&0.0151&0.0366&0.0333&0.0566&0.1599&0.0196&0.0337&0.0189&0.0770&0.0481\\
AD&0.0150&0.0352&0.0345&0.0568&0.1658&0.0200&0.0329&0.0188&0.0757&0.0459\\
\hline
&\multicolumn{5}{c}{\textbf{Local QMLE (}$\bolds{\bar{\theta}_n}$\textbf{)}}
&\multicolumn{5}{c}{\textbf{Local QMLE (}$\bolds{\bar{\theta}_n}$\textbf{)}}\\[-4pt]
& \multicolumn{5}{c}{\hrulefill} & \multicolumn{5}{c}{\hrulefill}\\
&\multicolumn{1}{c}{$\bolds{\bar{\mu}_n}$}&\multicolumn{1}{c}{$\bolds{\bar{\phi}_{1n}}$}&\multicolumn{1}{c}{$\bolds{\bar{\alpha}_{0n}}$}
&\multicolumn{1}{c}{$\bolds{\bar{\alpha}_{1n}}$}&\multicolumn{1}{c}{$\bolds{\bar{\beta}_{1n}}$}&\multicolumn{1}{c}{$\bolds{\bar{\mu}_{n}}$}
&\multicolumn{1}{c}{$\bolds{\bar{\phi}_{1n}}$}&\multicolumn{1}{c}{$\bolds{\bar{\alpha}_{0n}}$}&\multicolumn{1}{c}{$\bolds{\bar{\alpha}_{1n}}$}
&\multicolumn{1}{c@{}}{$\bolds{\bar{\beta}_{1n}}$}\\
\hline
Bias&0.0009&-0.0048&0.0055&0.0038&-0.0252&0.0004&-0.0031&0.0024&-0.0019&-0.0027\\
SD&0.0145&0.0300&0.0322&0.0454&0.1535&0.0195&0.0287&0.0183&0.0633&0.0442\\
AD&0.0145&0.0294&0.0320&0.0460&0.1517&0.0197&0.0279&0.0181&0.0644&0.0424\\
\hline
\end{tabular*}}}
\vspace*{-6pt}
\end{table}
%

%
%t3 ###
%ta3 #&#
\begin{table}
\tabcolsep=0pt
\caption{Estimators for model (\protect\ref{41}) when $\eta_t\sim t_3$}
\label{table3}
{\fontsize{8.5pt}{11pt}\selectfont{\begin{tabular*}{\tablewidth}
{@{\extracolsep{\fill}}ld{2.4}d{2.4}d{1.4}d{1.4}d{2.4}d{2.4}d{2.4}d{1.4}d{2.4}d{2.4}@{}}
\hline
&\multicolumn{5}{c}{$\bolds{\theta_0=(0.0, 0.5, 0.1, 0.18, 0.4)}$}
&\multicolumn{5}{c@{}}{$\bolds{\theta_0=(0.0, 0.5, 0.1, 0.2, 0.4)}$}\\[-4pt]
& \multicolumn{5}{c}{\hrulefill} & \multicolumn{5}{c}{\hrulefill}\\
&\multicolumn{5}{c}{\textbf{Self-weighted QMELE (}$\bolds{\hat{\theta
}_{sn}}$\textbf{)}} & \multicolumn{5}{c@{}}{\textbf{Self-weighted QMELE (}$\bolds{\hat{\theta
}_{sn}}$\textbf{)}}\\[-4pt]
& \multicolumn{5}{c}{\hrulefill} & \multicolumn{5}{c}{\hrulefill}\\
&\multicolumn{1}{c}{$\bolds{\hat{\mu}_{sn}}$}&\multicolumn{1}{c}{$\bolds{\hat{\phi}_{1sn}}$}&\multicolumn{1}{c}{$\bolds{\hat{\alpha}_{0sn}}$}
&\multicolumn{1}{c}{$\bolds{\hat{\alpha}_{1sn}}$}&\multicolumn{1}{c}{$\bolds{\hat{\beta}_{1sn}}$}&\multicolumn{1}{c}{$\bolds{\hat{\mu}_{sn}}$}
&\multicolumn{1}{c}{$\bolds{\hat{\phi}_{1sn}}$}&\multicolumn{1}{c}{$\bolds{\hat
{\alpha}_{0sn}}$} & \multicolumn{1}{c}{$\bolds{\hat{\alpha}_{1sn}}$}&\multicolumn{1}{c@{}}{$\bolds{\hat{\beta}_{1sn}}$}\\
\hline
Bias&0.0004&-0.0037&0.0059&0.0081&-0.0202&-0.0005&-0.0026&0.0032&0.0088&-0.0158\\
SD&0.0231&0.0416&0.0289&0.0600&0.1084&0.0221&0.0404&0.0252&0.0619&0.0968\\
AD&0.0233&0.0393&0.0282&0.0620&0.1101&0.0238&0.0393&0.0266&0.0637&0.1001\\
\hline
%&\multicolumn{5}{c|}{}&\multicolumn{5}{c}{}\\
&\multicolumn{5}{c}{\textbf{Local QMELE (}$\bolds{\hat{\theta}_n}$\textbf{)}}
&\multicolumn{5}{c}{\textbf{Local QMELE (}$\bolds{\hat{\theta}_n}$\textbf{)}}\\[-4pt]
& \multicolumn{5}{c}{\hrulefill} & \multicolumn{5}{c}{\hrulefill}\\
&\multicolumn{1}{c}{$\bolds{\hat{\mu}_n}$}&\multicolumn{1}{c}{$\bolds{\hat{\phi}_{1n}}$}&\multicolumn{1}{c}{$\bolds{\hat{\alpha}_{0n}}$}
&\multicolumn{1}{c}{$\bolds{\hat{\alpha}_{1n}}$}&\multicolumn{1}{c}{$\bolds{\hat{\beta}_{1n}}$}&\multicolumn{1}{c}{$\bolds{\hat{\mu}_{n}}$}
&\multicolumn{1}{c}{$\bolds{\hat{\phi}_{1n}}$}&\multicolumn{1}{c}{$\bolds{\hat{\alpha}_{0n}}$}&\multicolumn{1}{c}{$\bolds{\hat{\alpha}_{1n}}$}
&\multicolumn{1}{c@{}}{$\bolds{\hat{\beta}_{1n}}$}\\
\hline
Bias&0.0011&-0.0039&0.0041&0.0011&-0.0115&0.0001&-0.0028&0.0019&0.0029&-0.0092\\
SD&0.0229&0.0328&0.0256&0.0429&0.0955&0.0218&0.0325&0.0226&0.0450&0.0842\\
AD&0.0228&0.0314&0.0252&0.0461&0.0918&0.0233&0.0317&0.0243&0.0483&0.0851\\
\hline
&\multicolumn{5}{c}{\textbf{Self-weighted QMLE (}$\bolds{\bar{\theta
}_{sn}}$\textbf{)}}&\multicolumn{5}{c}{\textbf{Self-weighted QMLE
(}$\bolds{\bar{\theta}_{sn}}$\textbf{)}}\\[-4pt]
& \multicolumn{5}{c}{\hrulefill} & \multicolumn{5}{c}{\hrulefill}\\
&\multicolumn{1}{c}{$\bolds{\bar{\mu}_{sn}}$}&\multicolumn{1}{c}{$\bolds{\bar{\phi}_{1sn}}$}&\multicolumn{1}{c}{$\bolds{\bar{\alpha}_{0sn}}$}
&\multicolumn{1}{c}{$\bolds{\bar{\alpha}_{1sn}}$}&\multicolumn{1}{c}{$\bolds{\bar{\beta}_{1sn}}$}&\multicolumn{1}{c}{$\bolds{\bar{\mu}_{sn}}$}
&\multicolumn{1}{c}{$\bolds{\bar{\phi}_{1sn}}$}&\multicolumn{1}{c}{$\bolds{\bar{\alpha}_{0sn}}$}&\multicolumn{1}{c}{$\bolds{\bar{\alpha}_{1sn}}$}
&\multicolumn{1}{c@{}}{$\bolds{\bar{\beta}_{1sn}}$}\\
\hline
Bias&-0.0056&-0.0151&0.0029&0.0503&-0.0594&0.0036&-0.0141&0.0115&0.0442&-0.0543\\
SD&0.9657&0.1045&0.0868&0.2521&0.1740&0.1827&0.1065&0.3871&0.2164&0.1605\\
AD&0.0536&0.0907&33.031&0.1795&34.498&0.0517&0.0876&138.38&0.1875&11.302\\
\hline
&\multicolumn{5}{c}{\textbf{Local QMLE (}$\bolds{\bar{\theta}_n}$\textbf{)}}
&\multicolumn{5}{c}{\textbf{Local QMLE (}$\bolds{\bar{\theta}_n}$\textbf{)}}\\[-4pt]
& \multicolumn{5}{c}{\hrulefill} & \multicolumn{5}{c}{\hrulefill}\\
&\multicolumn{1}{c}{$\bolds{\bar{\mu}_n}$}&\multicolumn{1}{c}{$\bolds{\bar{\phi}_{1n}}$}&\multicolumn{1}{c}{$\bolds{\bar{\alpha}_{0n}}$}
&\multicolumn{1}{c}{$\bolds{\bar{\alpha}_{1n}}$}&\multicolumn{1}{c}{$\bolds{\bar{\beta}_{1n}}$}&\multicolumn{1}{c}{$\bolds{\bar{\mu}_{n}}$}
&\multicolumn{1}{c}{$\bolds{\bar{\phi}_{1n}}$}&\multicolumn{1}{c}{$\bolds{\bar{\alpha}_{0n}}$}&\multicolumn{1}{c}{$\bolds{\bar{\alpha}_{1n}}$}
&\multicolumn{1}{c@{}}{$\bolds{\bar{\beta}_{1n}}$}\\
\hline
Bias&-0.0048&-0.0216&-2.1342&0.0185&3.7712&-0.0010&-0.0203&1.3241&0.0253&-0.1333\\
SD&0.0517&0.1080&38.535&0.3596&83.704&0.0521&0.1318&42.250&0.2524&3.4539\\
AD&0.0508&0.0661&55.717&0.1447&45.055&0.0520&0.0707&13.761&0.1535&1.1343\\
\hline
\end{tabular*}}}
\end{table}

From Table \ref{table1}, when $\eta_t\sim \operatorname{Laplace}(0,1)$, we can see
that the
self-weighted QMELE has smaller AD and SD than those of both the
self-weighted QMLE and the local QMLE. When $\eta_t\sim N(0,1)$, in
Table \ref{table2}, we can see that the self-weighted QMLE has smaller
AD and
SD than those of both the self-weighted QMELE and the local QMELE.
From Table \ref{table3}, we note that the SD and AD of both the self-weighted
QMLE and the local QMLE are not close to each other since their
asymptotic variances are infinite, while the SD and AD of the
self-weighted QMELE and the local QMELE are very close to each
other. Except $\bar{\theta}_n$ in Table \ref{table3}, we can see that
all four
estimators in Tables \ref{table1}--\ref{table3} have very small biases,
and the local QMELE and local QMLE always have the smaller SD and AD than
those of the self-weighted QMELE and self-weighted QMLE,
respectively. This conclusion holds no matter with GARCH errors
(finite variance) or IGARCH errors. This coincides with what we
expected. Thus, if the tail index of the data is greater than 2 but
$E\eta_t^{4}=\infty$, we suggest to use the local QMELE in practice;
see also \citet{r25} for a~discussion.

Overall, the simulation results show that the self-weighted QMELE
and the local QMELE have a good performance in the finite sample,
especially for the heavy-tailed innovations.

%s5 ###
%se5 #&#
\section{A real example}\label{sec5}
\begin{figure}
\begin{tabular}{@{}cc@{}}

\includegraphics{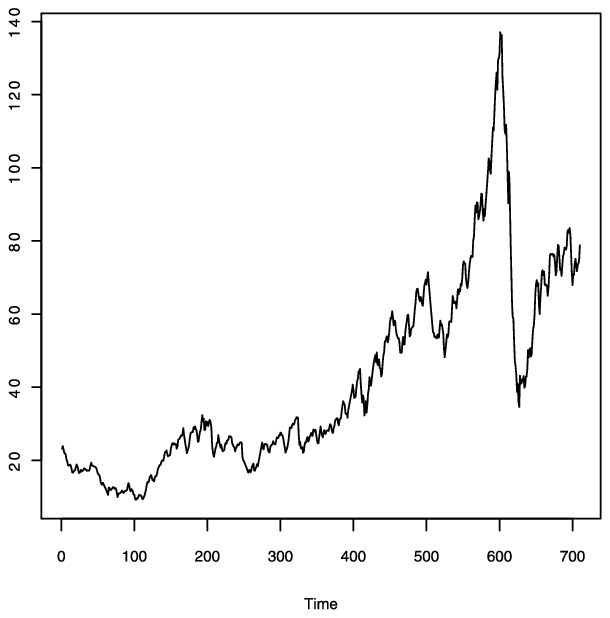}
 & \includegraphics{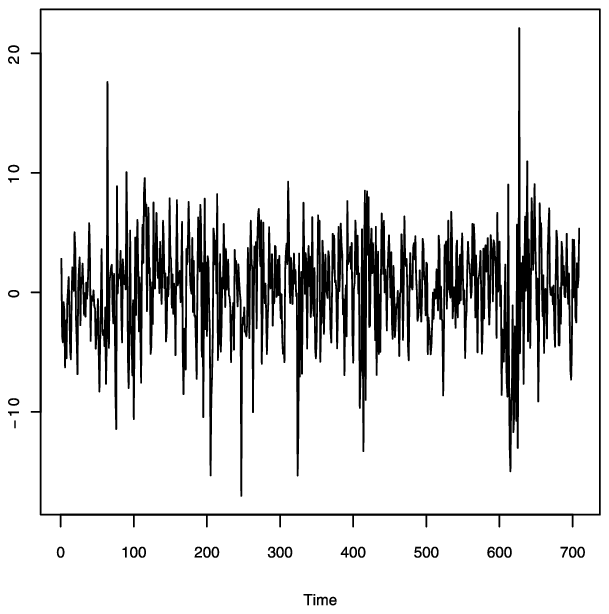}\\
(a) & (b)
\end{tabular}
\caption{\textup{(a)} The weekly world crude oil prices (dollars per barrel)
from January 3, 1997 to August~6, 2010 and \textup{(b)} its 100 times log
return.}\label{figure3}
\end{figure}

\begin{figure}
\begin{tabular}{@{}cc@{}}

\includegraphics{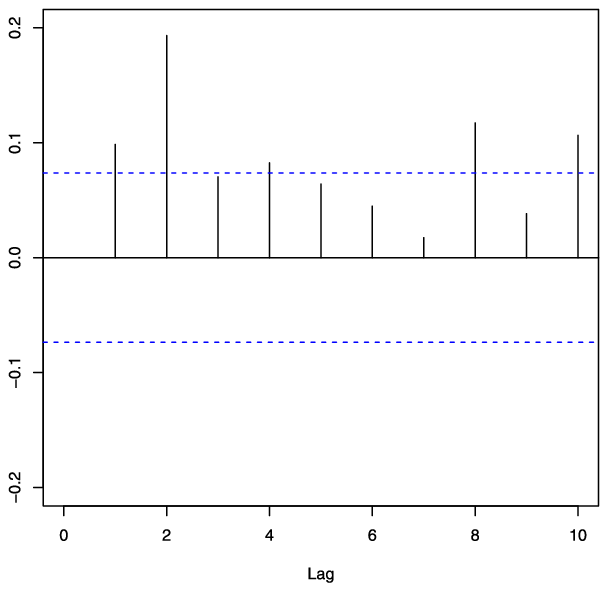}
 & \includegraphics{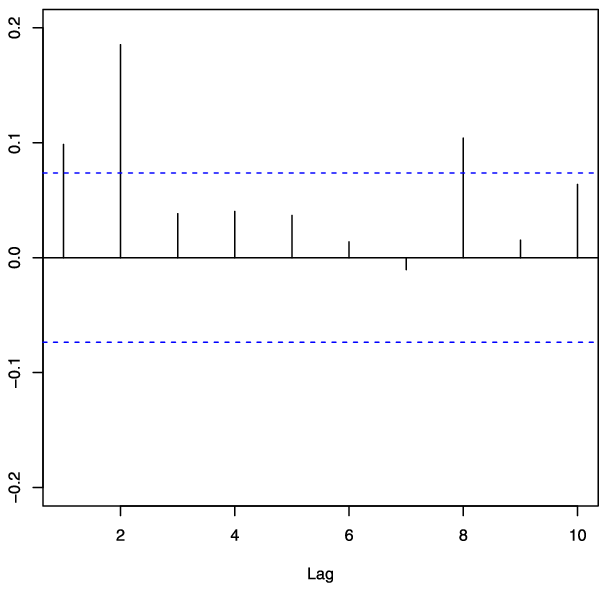}\\
(a) & (b)
\end{tabular}
\caption{\textup{(a)} The autocorrelations for $\{\hat{\e}_{t}^{2}\}$ and
\textup{(b)} the partial autocorrelations for $\{\hat{\e}_{t}^{2}\}$.}
\label{figure4}
\end{figure}

In this section, we study the weekly world crude oil price (dollars
per barrel) from January 3, 1997 to August 6, 2010, which has in
total 710 observations; see Figure \ref{figure3}(a). Its 100 times log-return,
denoted by $\{y_{t}\}_{t=1}^{709}$, is plotted in Figure \ref{figure3}(b). The
classic method based on the Akaike's information criterion (AIC)
leads to the following model:
%
%e5.1 ###
%e5.1 #&#
\begin{eqnarray} \label{51}
y_{t}&=& 0.2876\e_{t-1}+0.1524\e_{t-3}+\e_{t},\nonumber\\[-8pt]\\[-8pt]
&&(0.0357) \hspace*{24pt}(0.0357)\nonumber
\end{eqnarray}
where the standard errors are in parentheses, and the estimated
value of $\sigma^{2}_{\e}$ is 16.83. Model (\ref{51}) is stationary, and
none of the first ten autocorrelations or partial autocorrelations
of the residuals $\{\hat{\e}_{t}\}$ are significant at the 5\%
level. However, looking at the autocorrelations of
$\{\hat{\e}_{t}^{2}\}$, it turns out that the 1st, 2nd and 8th all
exceed two asymptotic standard errors; see Figure \ref{figure4}(a). Similar
results hold for the partial autocorrelations of
$\{\hat{\e}_{t}^{2}\}$ in Figure \ref{figure4}(b). This shows that
$\{\e_{t}^{2}\}$ may be highly correlated, and hence there may exist
ARCH effects.

%f3 ###
%fi3 #&#

%f4 ###
%fi4 #&#

\begin{figure}

\includegraphics{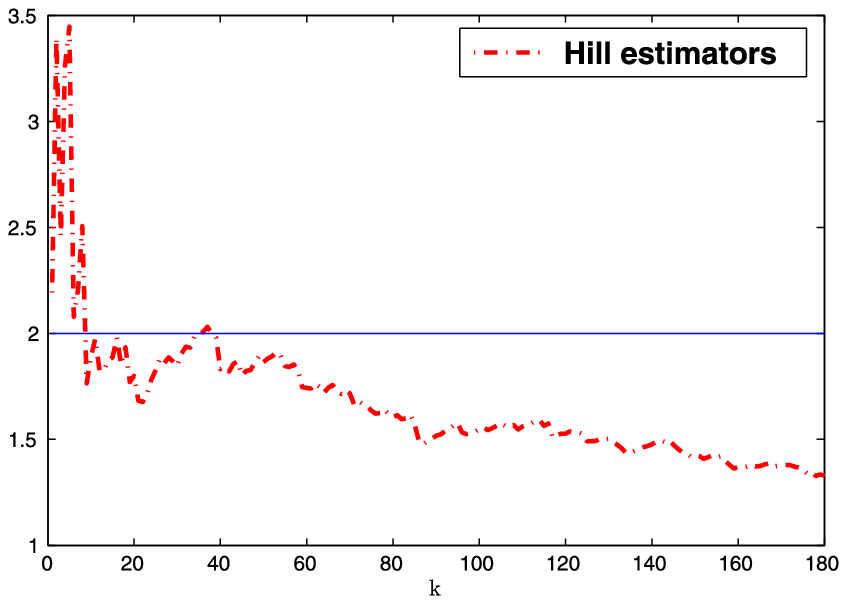}

\caption{Hill estimators $\{\hat{\alpha}_{y}(k)\}$ for
$\{y_{t}^{2}\}$.}
\label{figure5}
\end{figure}

\begin{figure}
\begin{tabular}{@{}cc@{}}

\includegraphics{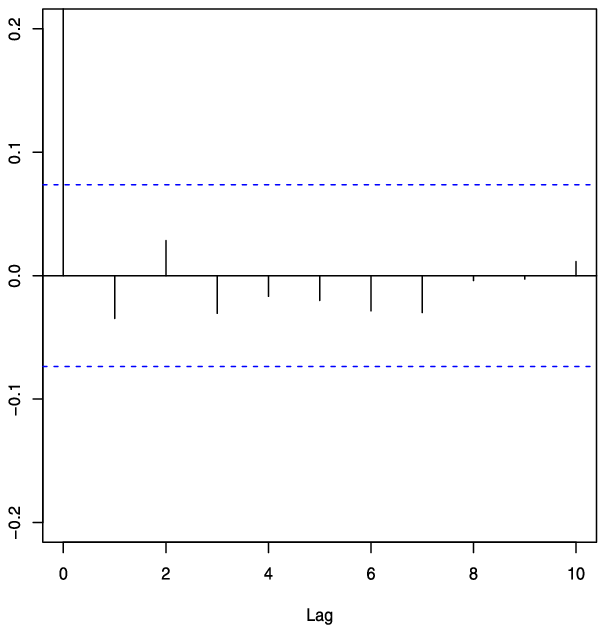}
 & \includegraphics{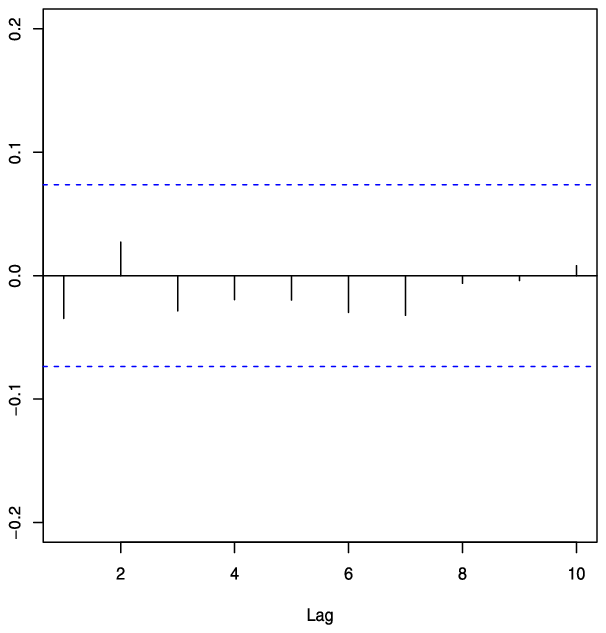}\\
(a) & (b)
\end{tabular}
\caption{\textup{(a)} The autocorrelations for $\{\hat{\eta}_{t}^{2}\}$ and
\textup{(b)} the partial autocorrelations for $\{\hat{\eta}_{t}^{2}\}$.}
\label{figure6}
\end{figure}

Thus, we try to use a MA(3)--GARCH$(1,1)$ model to fit the data set~$\{y_{t}\}$.
To begin with our estimation, we first estimate the
tail index of $\{y_{t}^{2}\}$ by using Hill's estimator
$\{\hat{\alpha}_{y}(k)\}$ with $k=1,\ldots,180$, based on
$\{y_{t}^{2}\}_{t=1}^{709}$. The plot of $\{\hat{\alpha}_{y}(k)\}_{k=1}^{180}$
is given in Figure \ref{figure5}, from which we can see that the tail
index of
$\{y_{t}^{2}\}$ is between 1 and 2, that is, $Ey_{t}^{4}=\infty$. So,
the standard QMLE procedure is not suitable. Therefore, we
first use the self-weighted QMELE to estimate the MA(3)--GARCH$(1,1)$
model, and then use the one-step iteration as in Section \ref{sec3} to obtain
its local QMELE. The fitted model is as follows:
%
%e5.2 ###
%e5.2 #&#
\begin{eqnarray}\label{52}
y_{t}&=&0.3276\e_{t-1}+0.1217\e_{t-3}+\e_{t},\nonumber\\[-2pt]
&&(0.0454)\hspace*{23.4pt} (0.0449)\nonumber\\[-10pt]\\[-10pt]
h_{t}&=&0.5147+0.0435\e_{t-1}^{2}+0.8756h_{t-1},\nonumber\\[-2pt]
&&(0.3248)\hspace*{4.2pt} (0.0159)\hspace*{24.5pt} (0.0530)\nonumber
\end{eqnarray}
where the standard errors are in parentheses. Again model (\ref{52}) is
stationary, and none of first ten autocorrelations or partial\vadjust{\eject}
autocorrelations of the residuals $\hat{\eta}_{t}\triangleq
\hat{\e}_{t}\hat{h}_{t}^{-1/2}$ are significant at the 5\% level.
Moreover, the first ten
autocorrelations and partial autocorrelations of $\{\hat{\eta
}_{t}^{2}\}$
are also within two asymptotic standard errors; see Figure \ref
{figure6}(a) and
(b). All these results suggest that model (\ref{52}) is adequate for the
data set $\{y_{t}\}$.

%f5 ###
%fi5 #&#

%f6 ###
%fi6 #&#

\begin{figure}

\includegraphics{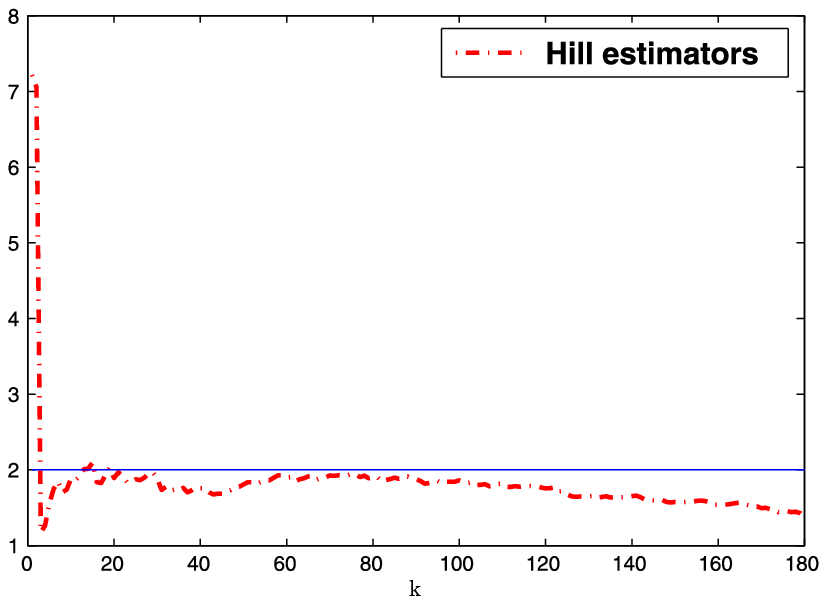}

\caption{The Hill estimators $\{\hat{\alpha}_{\eta}(k)\}$ for
$\{\hat{\eta}_{t}^{2}\}$ of model (\protect\ref{52}).}
\label{figure7}
\end{figure}

Finally, we estimate the tail index of $\eta_{t}^{2}$ in model (\ref{52})
by using Hill's estimator $\hat{\alpha}_{\eta}(k)$ with
$k=1,\ldots,180$, base on $\{\hat{\eta}_{t}^{2}\}$. The plot of
$\{\hat{\alpha}_{\eta}(k)\}_{k=1}^{180}$ is given in Figure \ref
{figure7}, from
%
%f7 ###
%fi7 #&#
which we can see that $E\eta_{t}^{2}$ is most likely finite, but~$E\eta_{t}^{4}$ is infinite.
Furthermore, the estimator of
$E\eta_{t}^{2}$ is $\sum_{t=1}^{n} \hat{\eta}_{t}^{2}/n=1.6994$, and
it turns out that $\hat{\alpha}_{1n}(\sum_{t=1}^{n}
\hat{\eta}_{t}^{2}/n)+\hat{\beta}_{1n}=0.9495$. This means
that \mbox{$E\e_{t}^{2}<\infty$}. Therefore, all the assumptions of Theorem
\ref{theorem31} are most likely satisfied. In particular, the estimated tail
indices of $\{y_{t}^{2}\}$ and $\{\hat{\eta}_{t}^{2}\}$ show the
evidence that the self-weighted/local QMELE is necessary in
modeling the crude oil price.

%s6 ###
%se6 #&#
\section{\texorpdfstring{Proofs of Lemmas \protect\ref{lemma22} and \protect\ref{lemma23}}{Proofs of Lemmas 2.2 and 2.3}}
\label{sec6}
In this section, we give the proofs of Lemmas \ref{lemma22} and \ref
{lemma23}. In the rest
of this paper, we denote $C$ as a universal constant, and $G(x)$ be the
distribution function of $\eta_t$.
\begin{pf*}{Proof of Lemma \ref{lemma22}}
A direct calculation gives
\[
\xi_t(u)=-u'\f{2w_t}{\sqrt{h_t}}\,\f{\p\e_t(\gamma_0)}{\p\theta}
M_t(u),
\]
where $M_t(u)=\int_{0}^1 X_t( -q_{1t}(u)s)\,ds$. Thus, we
have
\[
|\Pi_{1n}(u)|\leq2\|u\|\sum_{j=1}^m
\Biggl|\f{w_t}{\sqrt{h_t}}\,\f{\p\e_t(\gamma_0)}{\p\theta_j}\sum_{t=1}^n\{
M_t(u)-E[M_t(u)|\mathcal{F}_{t-1}]\}
\Biggr|.
\]
It is sufficient to show that
%
%e6.1 ###
%e6.1 #&#
\begin{equation}\label{61}\quad
\qquad \Biggl|\f{w_t}{\sqrt{h_t}}\,\f{\p\e_t(\theta_0)}{\p\theta_j}\sum_{t=1}^n\{
M_t(u_n)-E[M_t(u_n)|\mathcal{F}_{t-1}]\}\Biggr|=o_p\bigl(\sqrt{n}+n\|u_n\|\bigr),
\end{equation}
for each $1\leq j\leq m$.
Let $m_t=w_t h_t^{-1/2}\p\e_t(\theta_0)/\p\theta_j$, $
f_t(u)=m_tM_t(u)$ and
\[
D_n(u)=\f{1}{\sqrt{n}}\sum_{t=1}^n
\{f_t(u)-E[f_t(u)|\mathcal{F}_{t-1}]\}.
\]
Then, in order to prove (\ref{61}), we only need to show that for any
$\eta>0$,
%
%e6.2 ###
%e6.2 #&#
\begin{equation} \label{62}
\sup_{\|u\|\leq\eta}\f{|D_n(u)|}{1+\sqrt{n}\|u\|}=o_p(1).
\end{equation}
Note that $m_t=\max\{m_t, 0\}-\max\{-m_t, 0\}$. To make it simple,
we only prove the case when $m_t\ge0$.

We adopt the method in Lemma 4 of \citet{r32}. Let
$\mathfrak{F}=\{f_t(u): \|u\|\leq\eta\}$ be a collection of
functions indexed by $u$. We first verify that $\mathfrak{F}$
satisfies the bracketing condition in Pollard (\citeyear{r32}), page
304. Denote
$B_r(\zeta)$ be an open neighborhood of $\zeta$ with radius $r>0$.
For any fix $\e>0$ and $0<\delta\leq\eta$, there is a sequence of
small cubes $\{B_{\e\delta/C_1}(u_i)\}_{i=1}^{K_\e}$ to cover
$B_\delta(0)$, where $K_{\varepsilon}$ is an integer less than
$c_{0}\varepsilon^{-m}$ and $c_{0}$ is a constant not depending on
$\varepsilon$ and $\delta$; see Huber (\citeyear{r17}), page 227. Here,
$C_1$ is
a constant to be selected later. Moreover, we can choose
$U_i(\delta)\subseteq B_{\e\delta/C_1}(u_i)$ such that
$\{U_i(\delta)\}_{i=1}^{K_\e}$ be a partition of $B_\delta(0)$. For
each $u\in U_i(\delta)$, we define the bracketing functions as
follows:
\[
f^{\pm}_{t}(u)= m_t\int_{0}^1 X_t\biggl(-q_{1t}(u)s\pm\f{\e
\delta}{C_1\sqrt{h_t}} \biggl\|\f{\p\e_t(\gamma_0)}{\p\theta}\biggr\|
\biggr)\,ds.
\]
Since the indicator function is nondecreasing and $m_t\geq0$, we
can see that, for any $u\in U_i(\delta)$,
\[
f^{-}_{t}(u_i)\leq f_t(u)\leq f^{+}_{t}(u_i).
\]
Note that $\sup_{x\in R}g(x)<\infty$. It is straightforward to see
that
%
%e6.3 ###
%e6.3 #&#
\begin{equation}\label{63}
\qquad E[f^{+}_{t}(u_i)-f^{-}_{t}(u_i)|\mathcal{F}_{t-1}]\leq
\f{2\e\delta}{C_1}\sup_{x\in R}g(x)
\f{w_t}{h_t}\biggl\|\f{\p\e_t(\gamma_0)}{\p\theta}\biggr\|^2\equiv
\f{\e\delta\Delta_t}{C_1}.
\end{equation}
Setting $C_1=E(\Delta_t)$, we have
\[
E[f^{+}_{t}(u_i)-f^{-}_{t}(u_i)]=E\{
E[f^{+}_{t}(u_i)-f^{-}_{t}(u_i)|\mathcal{F}_{t-1}]\}\leq
\e\delta.
\]
Thus, the family $\mathfrak{F}$ satisfies
the bracketing condition.

Put $\delta_k=2^{-k}\eta$. Define $B(k)\equiv B_{\delta_k}(0)$, and
$A(k)$ to be the annulus $B(k)/\allowbreak B(k+1)$. Fix $\e>0$, for each $1\leq
i\leq K_\e$, by the bracketing condition, there exists a partition
$\{U_i(\delta_k)\}_{i=1}^{K_\e}$ of $B(k)$.

We first consider the upper tail. For $u\in U_i(\delta_k)$, by (\ref{63})
with $\delta=\delta_k$, we have
\begin{eqnarray*}
D_n(u)&\leq&\f{1}{\sqrt{n}}\sum_{t=1}^n \{
f^{+}_t(u_i)-E[f^{-}_t(u_i)|\mathcal{F}_{t-1}]\} \\
&=& D^{+}_n(u_i)+\f{1}{\sqrt{n}} \sum_{t=1}^n
E[f^{+}_t(u_i)-f^{-}_t(u_i)|\mathcal{F}_{t-1}]\\
&\leq& D^{+}_n(u_i)+\sqrt{n}\e\delta_k \Biggl[\f{1}{nC_1}\sum_{t=1}^n
\Delta_t\Biggr],
\end{eqnarray*}
where
\[
D^{+}_n(u_i)=\f{1}{\sqrt{n}}\sum_{t=1}^n \{
f^{+}_t(u_i)-E[f^{+}_t(u_i)|\mathcal{F}_{t-1}]\}.
\]
Denote the event
\[
E_n=\Biggl\{\omega: \f{1}{nC_1}\sum_{t=1}^n \Delta_t(\omega)<2 \Biggr\}.
\]
On $E_n$ with $u\in U_i(\delta_k)$, it follows that
%
%e6.4 ###
%e6.4 #&#
\begin{equation}\label{64}
D_n(u)\leq D^{+}_n(u_i)+2\sqrt{n}\e\delta_k.
\end{equation}
On $A(k)$, the divisor $1+\sqrt{n}\|u\|>
\sqrt{n}\delta_{k+1}=\sqrt{n}\delta_k/2$. Thus, by (\ref{64}) and
Chebyshev's inequality, it follows that
%
%e6.5 ###
%e6.5 #&#
\begin{eqnarray} \label{65}
&&P\biggl(\sup_{u\in A(k)}\f{D_n(u)}{1+\sqrt{n}\|u\|}>6\e, E_n\biggr) \nonumber\\
&&\qquad \leq P\Bigl(\sup_{u\in A(k)}D_n(u)>3\sqrt{n}\e\delta_k, E_n\Bigr) \nonumber\\
&&\qquad \leq P\Bigl(\max_{1\leq i\leq K_\e} \sup_{u\in
U_i(\delta_k)\cap A(k)} D_n(u)>3\sqrt{n}\e\delta_k, E_n\Bigr)
\nonumber\\[-8pt]\\[-8pt]
&&\qquad \leq P\Bigl(\max_{1\leq i\leq K_\e} D^{+}_n(u_i)>\sqrt{n}\e
\delta_k, E_n\Bigr) \nonumber\\
&&\qquad \leq K_\e\max_{1\leq i\leq K_\e}
P\bigl(D^{+}_n(u_i)>\sqrt{n}\e\delta_k\bigr) \nonumber\\
&&\qquad \leq K_\e\max_{1\leq i\leq K_\e}
\f{E[(D^{+}_n(u_i))^2]}{n\e^2\delta_k^2}.\nonumber
\end{eqnarray}
Note that $|q_{1t}(u_i)|\leq C\delta_k\xi_{\rho t-1}$ and $m^2_t\leq
Cw_t^2\xi_{\rho t-1}^2$ for some $\rho\in(0,1)$ by Lem\-ma~\ref{lemmaA1}(i), and $\sup_{x\in
R}g(x)<\infty$ by Assumption \ref{asm26}. By Taylor's expansion, we have
\begin{eqnarray*}
E[(f^{+}_t(u_i))^2]&=&E\{E[(f^{+}_t(u_i))^2|\mathcal{F}_{t-1}]\}\\
&\leq& E\biggl[m^2_t\int_{0}^1 E\biggl[\biggl|X_t\biggl(-q_{1t}(u_i)s+\f{\e\delta_k}{C_1\sqrt
{h_t}} \biggl\|\f{\p\e_t(\gamma_0)}{\p\theta}\biggr\| \biggr) \biggr|\Big|\mathcal{F}_{t-1}\biggr]\,ds \biggr] \\
&\leq& CE\Bigl[\sup_{|x|\leq\delta_k C\xi_{\rho t-1}}|G(x)-G(0)|w_t^2\xi
^2_{\rho t-1}\Bigr] \\
&\leq& \delta_kCE(w_t^2\xi_{\rho t-1}^{3}).
\end{eqnarray*}
Since $f^{+}_t(u_i)-E[f^{+}_t(u_i)|\mathcal{F}_{t-1}]$ is a
martingale difference sequence, by the previous inequality, it
follows that
%
%e6.6 ###
%e6.6 #&#
\begin{eqnarray}\label{66}
E[(D^{+}_n(u_i))^2]&=&\f{1}{n}\sum_{t=1}^n E\{
f^{+}_t(u_i)-E[f^{+}_t(u_i)|\mathcal{F}_{t-1}] \}^2 \nonumber\\
&\leq&\f{1}{n} \sum_{t=1}^n E[(f^{+}_t(u_i))^2] \nonumber\\[-8pt]\\[-8pt]
%&\leq\f{C}{n}\sum_{t=1}^n E[\sup_{|x|\leq\delta_k C'\xi_{\rho
%t-1}}|G(x)-G(0)|w_t^2\xi^2_{\rho t-1}],\nonumber\\
&\leq&\f{\delta_k}{n} \sum_{t=1}^n CE(w_t^2\xi_{\rho t-1}^{3}) \nonumber
\\
&\equiv&\pi_n(\delta_k).\nonumber
\end{eqnarray}
Thus, by (\ref{65}) and (\ref{66}), we have
\[
P\biggl(\sup_{u\in A(k)}\f{D_n(u)}{1+\sqrt{n}\|u\|}>6\e, E_n\biggr)\leq K_\e\f{\pi
_n(\delta_k)}{n\e^2\delta_k^2}.
\]
By a similar argument, we can get the same bound for the lower tail.
Thus, we can show that
%
%e6.7 ###
%e6.7 #&#
\begin{equation} \label{67}
P\biggl(\sup_{u\in A(k)}\f{|D_n(u)|}{1+\sqrt{n}\|u\|}>6\e,
E_n\biggr)\leq2K_\e\f{\pi_n(\delta_k)}{n\e^2\delta_k^2}.
\end{equation}

Since $\pi_n(\delta_k)\rightarrow0$ as $k\to\infty$, we can choose
$k_\e$ so that
\[
2\pi_n(\delta_k)K_\e/(\e\eta)^2<\e
\]
for $k\geq k_\e$. Let $k_n$ be an integer so that $n^{-1/2}\leq
2^{-k_n}< 2n^{-1/2}$. Split $\{u:\|u\|\leq\eta\}$ into two sets
$B(k_n+1)$ and $B(k_n+1)^c=\bigcup_{k=0}^{k_n} A(k)$.
By (\ref{67}), since $\pi_n(\delta_k)$ is bounded, we have
%
%e6.8 ###
%e6.8 #&#
\begin{eqnarray} \label{68}
&& P\biggl(\sup_{u\in
B(k_n+1)^c}\f{|D_n(u)|}{1+\sqrt{n}\|u\|}>6\e\biggr) \nonumber\\
&&\qquad \leq\sum_{k=0}^{k_n} P\biggl(\sup_{u\in A(k)}\f{|D_n(u)|}{1+\sqrt{n}\|u\|
}>6\e, E_n\biggr)+P(E_n^c) \\
&&\qquad \leq\f{1}{n}\sum_{k=0}^{k_\e-1} \f{C K_\e}{\e^2\eta^2}2^{2k}+\f{\e
}{n}\sum_{k=k_\e}^{k_n} 2^{2k}+P(E_n^c)\nonumber\\
&&\qquad \leq O\biggl(\f{1}{n}\biggr)+4\e\f{2^{2k_n}}{n}+P(E_n^c) \nonumber\\
&&\qquad \leq O\biggl(\f{1}{n}\biggr)+ 4\e+P(E_n^c).\nonumber
\end{eqnarray}

Since $1+\sqrt{n}\|u\|>1$ and $\sqrt{n}\delta_{k_n+1}<1$, using a
similar argument as for~(\ref{65}) together with (\ref{66}), we have
\begin{eqnarray*}
P\biggl(\sup_{u\in B(k_n+1)} \f{D_n(u)}{1+\sqrt{n}\|u\|}>3\e, E_n\biggr) &\leq&
P\Bigl(\max_{1\leq i\leq K_\e} D^{+}_n(u_i)>\e, E_n\Bigr) \\
&\leq&\f{K_\e\pi_n(\delta_{k_n+1})}{\e^2}.
\end{eqnarray*}
We can get the same bound for the lower tail. Thus, we have
%
%e6.9 ###
%e6.9 #&#
\begin{eqnarray}\label{69}
&& P\biggl(\sup_{u\in B(k_n+1)} \f{|D_n(u)|}{1+\sqrt{n}\|u\|}>3\e\biggr) \nonumber\\
&&\qquad =P\biggl(\sup_{u\in B(k_n+1)} \f{|D_n(u)|}{1+\sqrt{n}\|u\|}>3\e,
E_n\biggr)+P(E_n^c)\\
&&\qquad \le\f{2K_\e\pi_n(\delta_{k_n+1})}{\e^2}+P(E_n^c).\nonumber
\end{eqnarray}
Note that $\pi_n(\delta_{k_n+1})\rightarrow0$ as $n\to\infty$.
Furthermore, $P(E_n)\rightarrow1$ by the ergodic theorem. Hence,
\[
P(E_n^c)\to0 \qquad\mbox{as } n\to\infty.
\]
Finally, (\ref{62}) follows by (\ref{68}) and (\ref{69}). This completes
the proof.
\end{pf*}
\begin{pf*}{Proof of Lemma \ref{lemma23}}
(i). By a direct calculation, we
have
%
%e6.10 ###
%e6.10 #&#
\begin{eqnarray}\label{610}
\Pi_{2n}(u)&=&2\sum_{t=1}^n w_t\int_0^{-q_{1t}(u)} G(s)-G(0) \,ds \nonumber\\
&=& 2\sum_{t=1}^n w_t\int_0^{-q_{1t}(u)} sg(\varsigma^*) \,ds \\
&=& \bigl(\sqrt{n}u\bigr)' [K_{1n}+K_{2n}(u)]\bigl(\sqrt{n}u\bigr),\nonumber
\end{eqnarray}
where $\varsigma^*$ lies between 0 and $s$, and
\begin{eqnarray*}
K_{1n}&=&\f{g(0)}{n} \sum_{t=1}^n \f{w_t}{h_t(\theta_0)} \,\f{\p\e_t(\gamma
_0)}{\p\theta} \,\f{\p\e_t(\gamma_0)}{\p\theta'},\\
K_{2n}(u)&=&\f{2}{n\|u\|^2} \sum_{t=1}^n w_t\int_0^{-q_{1t}(u)}
s[g(\varsigma^*)-g(0)]\,ds.
\end{eqnarray*}
By the ergodic theorem, it is easy to see that
%
%e6.11 ###
%e6.11 #&#
\begin{equation}\label{611}
K_{1n}=\Sigma_1+o_p(1).
\end{equation}
Furthermore, since $|q_{1t}(u)|\leq C\|u\|\xi_{\rho t-1}$ for some
$\rho\in(0,1)$ by Lemma \ref{lemmaA1}(i), it is straightforward to see that
for any $\eta>0$,
\begin{eqnarray*}
\sup_{\|u\|\leq\eta}|K_{2n}(u)|&\leq&\sup_{\|u\|\leq
\eta}\f{2}{n\|u\|^2} \sum_{t=1}^n
w_t\int_{-|q_{1t}(u)|}^{|q_{1t}(u)|}
s|g(\varsigma^*)-g(0)|\,ds\\
&\leq&\f{1}{n} \sum_{t=1}^n \Bigl[\sup_{|s|\leq C\eta\xi_{\rho t-1}}
|g(s)-g(0)| w_t\xi_{\rho t-1}^2\Bigr].
\end{eqnarray*}
By Assumptions \ref{asm24} and \ref{asm26}, $E(w_t\xi_{\rho
t-1}^2)<\infty$ and $\sup_{x\in R} g(x)<\infty$. Then, by the dominated
convergence
theorem, we have
\[
\lim_{\eta\to0}E\Bigl[\sup_{|s|\leq C\eta\xi_{\rho t-1}}
|g(s)-g(0)|w_t\xi_{\rho t-1}^2 \Bigr]=0.
\]
Thus, by the stationarity of $\{y_{t}\}$ and Markov's theorem, for
$\forall\e,\delta>0$, $\exists\eta_{0}(\e)>0$, such that
%
%e6.12 ###
%e6.12 #&#
\begin{equation}\label{612}
P\Bigl(\sup_{\|u\|\leq
\eta_{0}}|K_{2n}(u)|>\delta\Bigr)<\frac{\e}{2}
\end{equation}
for all $n\geq1$. On the other hand, since $u_{n}=o_{p}(1)$, it
follows that
%
%e6.13 ###
%e6.13 #&#
\begin{equation}\label{613}
P(\|u_{n}\|>\eta_{0})<\frac{\e}{2}
\end{equation}
as $n$ is large enough. By (\ref{612}) and (\ref{613}), for $\forall\e
,\delta>0$, we have
\begin{eqnarray*}
P\bigl(|K_{2n}(u_{n})|>\delta\bigr)&\leq&
P\bigl(|K_{2n}(u_{n})|>\delta, \|u_{n}\|\leq
\eta_{0}\bigr)+
P(\|u_{n}\|>\eta_{0}) \\
&<& P\Bigl(\sup_{\|u\|\leq\eta_{0}}|K_{2n}(u)|>\delta\Bigr)+\frac{\e}{2}\\
&<&\e
\end{eqnarray*}
as $n$ is large enough, that is, $K_{2n}(u_{n})=o_{p}(1)$. Furthermore,
combining (\ref{610}) and (\ref{611}), we can see that (i) holds.

(ii) Let $\Pi_{3n}(u)=(\sqrt{n}u)'K_{3n}(\xi^{*})(\sqrt
{n}u)+K_{4n}(u)$, where
\begin{eqnarray*}
K_{3n}(\xi^{*})&=&\f{1}{n}\sum_{t=1}^n \f{w_t}{\sqrt{h_t}}\,\f{\p^2\e_t(\xi
^*)}{\p\theta\,\p\theta'}[I(\eta_t>0)-I(\eta_t<0)],\\
K_{4n}(u)&=&2\sum_{t=1}^n w_t\int_{-q_{1t}(u)}^{-q_t(u)} X_t(s) \,ds.
\end{eqnarray*}
By Assumption \ref{asm24} and Lemma \ref{lemmaA1}(i), there exists a constant
$\rho\in(0,1)$ such that
\[
E\biggl(\sup_{\xi^{*}\in\Lambda} \f{w_t}{\sqrt{h_t}}\biggl|\f{\p^2\e_t(\xi^*)}{\p
\theta\,\p\theta'}
[I(\eta_t>0)-I(\eta_t<0)]\biggr|\biggr)\leq CE(w_t\xi_{\rho t-1})<\infty.
\]
Since $\eta_{t}$ has median 0, the conditional expectation property
gives
\[
E\biggl(\f{w_t}{\sqrt{h_t}}\,\f{\p^2\e_t(\xi^*)}{\p\theta\,\p\theta'}[I(\eta
_t>0)-I(\eta_t<0)]\biggr)=0.
\]
Then, by Theorem 3.1 in \citet{r28}, we have
\[
\sup_{\xi^{*}\in\Lambda} |K_{3n}(\xi^{*})|=o_p(1).
\]
On the other hand,
\begin{eqnarray*}
\f{K_{4n}(u)}{n\|u\|^2}&=&\f{2}{n}\sum_{t=1}^n w_t
\int_0^{-q_{2t}(u)/\|u\|^2} X_t\bigl(\|u\|^2s-q_{1t}(u)\bigr)\,ds\\
&\equiv&\f{2}{n}\sum_{t=1}^{n} J_{1t}(u).
\end{eqnarray*}
By Lemma \ref{lemmaA1}, we have $|\|u\|^{-2}q_{2t}(u)|\leq C\xi_{\rho t-1}$
and $|q_{1t}(u)|\leq C\|u\|\xi_{\rho t-1}$ for some $\rho\in(0,1)$.
Then, for any $\eta>0$, we have
\begin{eqnarray*}
\sup_{\|u\|\leq\eta}|J_{1t}(u)| &\leq& w_t \int_{-C\xi_{\rho
t-1}}^{C\xi_{\rho t-1}} \{X_t(C\eta^2\xi_{\rho
t-1}+C\eta\xi_{\rho t-1})
\\
&&\hspace*{52.4pt}{} -X_t(-C\eta^2\xi_{\rho t-1}-C\eta\xi_{\rho
t-1})\}\,ds
\\
&\leq& 2Cw_t\xi_{\rho t-1} \{X_t(C\eta^2\xi_{\rho
t-1}+C\eta\xi_{\rho t-1})
\\
&&\hspace*{52.4pt}{} -X_t(-C\eta^2\xi_{\rho t-1}-C\eta\xi_{\rho
t-1})\}.
\end{eqnarray*}
By Assumptions \ref{asm24} and \ref{asm26} and the double expectation
property, it
follows that
\begin{eqnarray*}
E\Bigl[\sup_{\|u\|\leq\eta}|J_{1t}(u)|\Bigr]&\leq&
2CE[w_t\xi_{\rho t-1} \{G(C\eta^2\xi_{\rho
t-1}+C\eta\xi_{\rho t-1})
\\
&&\hspace*{65pt}{} -G(-C\eta^2\xi_{\rho
t-1}-C\eta\xi_{\rho t-1})\}] \\
&\leq& C(\eta^2+\eta)\sup_{x}g(x)E(w_t\xi_{\rho t-1}^{2})\to0
\end{eqnarray*}
as $\eta\to0$. Thus, as for (\ref{612}) and (\ref{613}), we can show that
$K_{4n}(u_{n})=o_{p}(n\|u_{n}\|^{2})$. This completes the proof
of\vspace*{2pt}
(ii).

(iii) Let $\Pi_{4n}(u)=(\sqrt{n}u)'[n^{-1}\sum_{t=1}^n
J_{2t}(\zeta^{*})](\sqrt{n}u)$, where
\[
J_{2t}(\zeta^*)=w_t
\biggl(\f{3}{8}\biggl|\f{\e_t(\gamma_0)}{\sqrt{h_t(\zeta^*)}}\biggr|
-\f{1}{4}\biggr)\f{1}{h_t^2(\zeta^*)} \,\f{\p h_t(\zeta^*)}{\p\theta}\,
\f{\p h_t(\zeta^*)}{\p\theta'}.
\]
By Assumption \ref{asm24} and Lemma \ref{lemmaA1}(ii)--(iv), there
exists a constant
$\rho\in(0,1)$ and a neighborhood $\Theta_{0}$ of $\theta_{0}$ such
that
\[
E\Bigl[\sup_{\zeta^{*}\in\Theta_{0}} |J_{2t}(\zeta^*)|
\Bigr]\leq CE[w_t\xi_{\rho t-1}^2(|\eta_t|\xi_{\rho
t-1}+1)]<\infty.
\]
Then, by Theorem 3.1 of \citet{r28}, we have
\[
\sup_{\zeta^{*}\in\Theta_{0}}\Biggl|\f{1}{n}\sum_{t=1}^n
J_{2t}(\zeta^{*})-E[J_{2t}(\zeta^{*})]\Biggr|=o_p(1).
\]
Moreover, since $\zeta^{*}_{n}\to\theta_{0}$ a.s., by the dominated
convergence theorem, we have
\[
\lim_{n\to\infty} E[J_{2t}(\zeta^{*}_{n})]=E[J_{2t}(\theta_0)]=\Sigma_{2}.
\]
Thus, (iii) follows from the previous two equations. This completes
the proof of~(iii).

(iv) Since $E|\eta_t|=1$, a similar argument as for part (iii)
shows that (iv) holds.

\mbox{}\phantom{i}(v) By Taylor's expansion, we have
\[
\f{1}{\sqrt{h_t(\theta_0+u)}}-\f{1}{\sqrt{h_t(\theta_0)}}=\f
{-u'}{2(h_t(\zeta^{*}))^{3/2}}\,\f{\p
h_t(\zeta^*)}{\p\theta},
\]
where $\zeta^*$ lies between $\theta_0$ and $\theta_0+u$. By
identity (\ref{26}), it is easy to see that
\begin{eqnarray*}
|\e_t(\gamma_0+u_1)|-|\e_t(\gamma_0)|&=&u'\,\f{\p\e_t(\xi^*)}{\p\theta
}[I(\eta_t>0)-I(\eta_t<0)]\\
&&{} +2u'\,\f{\p\e_t(\xi^*)}{\p\theta}\int_{0}^{1}
X_t\biggl(-\f{u'}{\sqrt{h_t}}\,\f{\p\e_t(\xi^*)}{\p\theta}s\biggr)\,ds,
\end{eqnarray*}
where $\xi^*$ lies between $\gamma_0$ and $\gamma_0+u_1$. By the
previous two equations, it follows that
\[
\sum_{t=1}^n
w_tC_t(u)=\bigl(\sqrt{n}u\bigr)'[K_{5n}(u)+K_{6n}(u)]\bigl(\sqrt{n}u\bigr),
\]
where
\begin{eqnarray*}
K_{5n}(u)&=&\f{1}{n}\sum_{t=1}^n \f{w_t}{2h_t^{3/2}(\zeta^{*})}\,\f{\p
h_t(\zeta^*)}{\p\theta} \,\f{\p\e_t(\xi^*)}{\p\theta'}
[I(\eta_t<0)-I(\eta_t>0)],\\
K_{6n}(u)&=&-\f{1}{n}\sum_{t=1}^n \f{w_t}{h_t^{3/2}(\zeta^{*})}\,\f{\p
h_t(\zeta^*)}{\p\theta} \,\f{\p\e_t(\xi^*)}{\p\theta'}\int_{0}^{1}
X_t\biggl(-\f{u'}{\sqrt{h_t}}\,\f{\p\e_t(\xi^*)}{\p\theta}s\biggr)\,ds.
\end{eqnarray*}
By Lemma \ref{lemmaA1}(i), (iii), (iv) and a similar argument as for part
(ii), it is easy to see that $K_{5n}(u_n)=o_p(1)$ and
$K_{6n}(u_n)=o_p(1)$. Thus, it follows that (v) holds. This
completes all of the proofs.
\end{pf*}

%s7 ###
%se7 #&#
\section{Concluding remarks}\label{sec7}
In this paper, we first propose a self-weighted QMELE for the
ARMA--GARCH model. The strong consistency and asymptotic normality of
the global self-weighted QMELE\vspace*{1pt} are established under a~fractional
moment condition of $\varepsilon_{t}$ with $E\eta_{t}^{2}<\infty$.
Based on this estimator, the local QMELE is showed to be asymptotically
normal for the ARMA--GARCH (finite variance) and --IGARCH models. The
empirical study shows that the self-weighted/local QMELE has a better
performance than the self-weighted/local QMLE when $\eta_{t}$ has a
heavy-tailed distribution, while the local QMELE is more efficient than
the self-weighted QMELE for the cases with a finite variance and
--IGARCH errors. We also give a real example to illustrate that our new
estimation procedure is necessary. According to our limit experience,
the estimated tail index of most of data sets lies in $[2,4)$ in
economics and finance. Thus, the local QMELE may be the most suitable
in practice if there is a further evidence to show that
$E\eta_{t}^{4}=\infty$.

%apA #&#
\begin{appendix}\label{app}
\section*{Appendix}

The Lemma \ref{lemmaA1} below is from \citet{r25}.
%
%leA.1 #&#
\begin{lem}\label{lemmaA1}
Let $\xi_{\rho t}$ be defined as in Assumption \ref{asm24}.
If Assumptions~\ref{asm21} and~\ref{asm22} hold, then there exists a constant $\rho\in (0,1)$ and a
neighborhood~$\Theta_0$ of $\theta_0$ such that:
\begin{eqnarray*}
\sup_{\Theta}|\e_{t-1}(\gamma)|&\leq& C\xi_{\rho t-1},\\
\mbox{\textup{(i)}\hspace*{37.5pt}\quad}\sup_{\Theta}
\biggl\|\f{\p\e_t(\gamma)}{\p\gamma}\biggr\|&\leq& C\xi_{\rho t-1}\quad \mbox{and}\\
\sup_{\Theta}\biggl\|\f{\p^2\e_t(\gamma)}{\p\gamma\,\p\gamma'}\biggr\|&\leq&
C\xi_{\rho t-1},\\
\mbox{\textup{(ii)}\hspace*{55pt}\quad} \sup_{\Theta} h_t(\theta)&\leq& C\xi^2_{\rho t-1}, \\
\mbox{\textup{(iii)}\quad\hspace*{11pt}} \sup_{\Theta_0}\biggl\|\f{1}{h_t(\theta)}\,\f{\p
h_t(\theta)}{\p\delta}\biggr\|&\leq& C\xi^{\iota_1}_{\rho
t-1}\qquad\mbox{for any } \iota_1\in(0,1),\\
\mbox{\textup{(iv)}\quad} \sup_{\Theta}\biggl\|\f{1}{\sqrt{h_t(\theta)}}\,\f{\p
h_t(\theta)}{\p\gamma}\biggr\|&\leq& C\xi_{\rho t-1}.
\end{eqnarray*}
\end{lem}
%
%leA.2 #&#
\begin{lem}\label{lemmaA2}
For any $\theta^*\in\Theta$, let
$B_\eta(\theta^*)=\{\theta\in\Theta:\|\theta-\theta^*\|<\eta\}$ be
an open neighborhood of $\theta^*$ with radius $\eta>0$. If
Assumptions \ref{asm21}--\ref{asm25} hold, then:
\begin{eqnarray*}
&&\mbox{\hphantom{ii}\textup{(i)}\quad} E\Bigl[\sup_{\theta\in\Theta} w_tl_t(\theta)\Bigr]<\infty,\\
&&\mbox{\hphantom{i}\textup{(ii)}\quad} E[w_tl_t(\theta)] \qquad\mbox{has a unique minimum at } \theta_0,\\
&&\mbox{\textup{(iii)}\quad} E\Bigl[\sup_{\theta\in
B_\eta(\theta^*)}w_t|l_t(\theta)-l_t(\theta^*)|\Bigr]\to0 \qquad\mbox{as } \eta\to0.
\end{eqnarray*}
\end{lem}
\begin{pf}
First, by (A.13) and (A.14) in \citet{r25} and Assumptions~\ref{asm24} and~\ref{asm25}, it follows that
\[
E\biggl[\sup_{\theta\in\Theta}
\f{w_t|\e_t(\gamma)|}{\sqrt{h_t(\theta)}}\biggr]\leq
CE[w_t\xi_{\rho t-1}(1+|\eta_t|)]<\infty
\]
for some
$\rho\in(0,1)$, and
\[
E\Bigl[\sup_{\theta\in\Theta}w_t\log
\sqrt{h_t(\theta)}\Bigr]<\infty;
\]
see Ling (\citeyear{r25}), page 864. Thus, (i)
holds.

Next, by a direct calculation, we have
\begin{eqnarray*}
E[w_tl_t(\theta)]&=&E\biggl[w_t\log\sqrt{h_t(\theta)}+\frac{w_t|\e_t(\gamma
_0)+(\gamma-\gamma_0)'({\p\e_t(\xi^*)}/{\p\theta})|}{\sqrt{h_t(\theta
)}}\biggr]\\
&=&E\biggl[w_t\log\sqrt{h_t(\theta)}+\frac{w_t}{\sqrt{h_t(\theta)}}E\biggl\{\!\biggl|\e
_t(\gamma_0)+(\gamma-\gamma_0)'\,\f{\p\e_t(\xi^*)}{\p\theta}\biggr|\Big|\mathcal
{F}_{t-1}\!\biggr\}\!\biggr]\\
&\geq& E\biggl[w_t\log\!\sqrt{h_t(\theta)}+\frac{w_t}{\sqrt{h_t(\theta)}}E(|\e
_t||\mathcal{F}_{t-1})\biggr]\\
&=&E\Biggl[w_t\Biggl(\log\sqrt{\f{h_t(\theta)}{h_t(\theta_0)}}+\sqrt{\f{h_t(\theta
_0)}{h_t(\theta)}}\Biggr)\Biggr]+E
\bigl[w_t\log\sqrt{h_t(\theta_0)}\bigr],
\end{eqnarray*}
where the last inequality holds since $\eta_t$ has a unique median
0, and obtains the minimum if and only if $\gamma=\gamma_0$ a.s.;
see \citet{r25}. Here, $\xi^*$ lies between~$\gamma$ and $\gamma_0$.
Considering the function $f(x)=\log x+a/x$ when $a\geq0$, it
reaches the minimum at $x=a$. Thus, $E[w_tl_t(\theta)]$ reaches the
minimum if and only if $\sqrt{h_t(\theta)}=\sqrt{h_t(\theta_0})$
a.s., and hence $\theta=\theta_0$; see \citet{r25}. Thus, we can
claim that $E[w_tl_t(\theta)]$ is uniformly minimized at $\theta_0$,
that is, (ii) holds.

Third, let $\theta^*=({\gamma^{*}}', {\delta^{*}}')'\in\Theta$. For
any $\theta\in B_\eta(\theta^*)$, using Taylor's expansion, we can
see that
\[
\log\sqrt{h_t(\theta)}-\log\sqrt{h_t(\theta^*)}=\f{(\theta-\theta
^*)'}{2h_t(\theta^{**})}\,\f{\p
h_t(\theta^{**})}{\p\theta},
\]
where $\theta^{**}$ lies between $\theta$ and $\theta^{*}$. By Lemma
\ref{lemmaA1}(iii)--(iv) and Assumption~\ref{asm24}, for some $\rho\in
(0,1)$, we have
\[
E\Bigl[\sup_{\theta\in B_{\eta}(\theta^*)}
w_t\bigl|\log\sqrt{h_t(\theta)}-\log\sqrt{h_t(\theta^*)}\bigr|\Bigr]\leq
C\eta E(w_t\xi_{\rho t-1})\rightarrow0
\]
as $\eta\rightarrow0$. Similarly,
\begin{eqnarray*}
E\biggl[\sup_{\theta\in B_{\eta}(\theta^*)}\f{w_t}{\sqrt{h_t(\theta)}}\bigl||\e
_t(\gamma)|-|\e_t(\gamma^*)|\bigr|\biggr]&\rightarrow&0 \qquad\mbox{as } \eta\rightarrow
0,\\
E\biggl[\sup_{\theta\in B_{\eta}(\theta^*)}
w_t|\e_t(\gamma^*)|\biggl|\f{1}{\sqrt{h_t(\theta)}}-\f{1}{\sqrt{h_t(\theta
^*)}}\biggr|\biggr]&\rightarrow&
0\qquad \mbox{as } \eta\rightarrow0.
\end{eqnarray*}
Then, it follows that (iii) holds. This completes all of the proofs
of Lemma~\ref{lemmaA2}.
\end{pf}
\begin{pf*}{Proof of Theorem \ref{theorem21}}
We use the method in \citet{r17}.
Let $V$ be any open neighborhood of $\theta_0\in\Theta$. By Lemma
\ref{lemmaA2}(iii), for any $\theta^*\in V^c=\Theta/V$ and $\e>0$, there
exists an $\eta_0>0$ such that
%
%e7.1 ###
%eA.1 #&#
\begin{equation}\label{A1}
E\Bigl[\inf_{\theta\in B_{\eta_0}(\theta^*)} w_tl_t(\theta)\Bigr]
\geq E[w_tl_t(\theta^*)]-\e.
\end{equation}
From Lemma \ref{lemmaA2}(i), by the ergodic theorem, it follows that
%
%e7.2 ###
%eA.2 #&#
\begin{equation}\label{A2}
\f{1}{n}\sum_{t=1}^n \inf_{\theta\in B_{\eta_0}(\theta^*)}
w_tl_t(\theta) \geq E\Bigl[\inf_{\theta\in B_{\eta_0}(\theta^*)}
w_t l_t(\theta)\Bigr]-\e
\end{equation}
as $n$ is large enough. Since $V^c$ is compact, we can choose
$\{B_{\eta_{0}}(\theta_i)\dvtx\theta_i\in V^c, i=1,2,\ldots, k\}$ to be a
finite covering of $V^c$. Thus, from (\ref{A1}) and (\ref{A2}), we have
%
%e7.3 ###
%eA.3 #&#
\begin{eqnarray}\label{A3}
\inf_{\theta\in V^c} L_{sn}(\theta)&=&\min_{1\leq i\leq k}\inf_{\theta
\in B_{\eta_0}(\theta_i)} L_{sn}(\theta) \nonumber\\
&\geq&\min_{1\leq i\leq k} \f{1}{n} \sum_{t=1}^n \inf_{\theta\in B_{\eta
_0}(\theta_i)} w_t l_t(\theta) \\
&\geq&\min_{1\leq i\leq k} E\Bigl[\inf_{\theta\in
B_{\eta_0}(\theta_i)} w_t l_t(\theta)\Bigr]-\e\nonumber
\end{eqnarray}
as $n$ is large enough. Note that the infimum on the compact set
$V^c$ is attained. For each $\theta_i\in V^c$, from Lemma \ref{lemmaA2}(ii),
there exists an $\e_0>0$ such that
%
%e7.4 ###
%eA.4 #&#
\begin{equation}\label{A4}
E\Bigl[\inf_{\theta\in B_{\eta_0}(\theta_i)} w_t
l_t(\theta)\Bigr]\geq E[w_t l_t(\theta_0)]+3\e_0.
\end{equation}
Thus, from (\ref{A3}) and (\ref{A4}), taking $\e=\e_0$, it follows that
%
%e7.5 ###
%eA.5 #&#
\begin{equation}\label{A5}
\inf_{\theta\in V^c} L_{sn}(\theta)\geq E[w_t l_t(\theta_0)]+2\e_0.
\end{equation}
On the other hand, by the ergodic theorem, it follows that
%
%e7.6 ###
%eA.6 #&#
\begin{equation} \label{A6}
\inf_{\theta\in V} L_{sn}(\theta)\leq
L_{sn}(\theta_0)=\f{1}{n}\sum_{t=1}^n w_t l_t(\theta_0)\leq E[w_t
l_t(\theta_0)]+\e_0.
\end{equation}
Hence, combining (\ref{A5}) and (\ref{A6}), it gives us
\[
\inf_{\theta\in V^c} L_{sn}(\theta)\geq E[w_t
l_t(\theta_0)]+2\e_0>E[w_t l_t(\theta_0)]+\e_0\geq\inf_{\theta\in
V} L_{sn}(\theta),
\]
which implies that
\[
\hat{\theta}_{sn}\in V\qquad\mbox{ a.s. for }\forall V\mbox{,  as } n \mbox
{ is large enough.}
\]
By the arbitrariness of $V$, it yields $\hat{\theta}_{sn}\rightarrow
\theta_0$ a.s. This completes the proof.
\end{pf*}
\end{appendix}

\section*{Acknowledgments}
The authors greatly appreciate the very helpful comments of two
anonymous referees, the Associate Editor and the Editor T.~Tony, Cai.

%suskaldyti doi

% imsref loaded by lrinkeviciute, 2011-07-27 08:45:58
% imsref loaded by lrinkeviciute, 2011-07-27 08:49:57
%

%
\printaddresses

\end{document}